\documentclass[11pt,twoside]{article}

\usepackage{amsbsy,amsfonts,amsmath,amssymb,eucal,mathrsfs}
\usepackage[all]{xy}
\usepackage[T1]{fontenc}

\def\itemn#1{\item[\hspace{0.6mm} {\rm (#1)}]}

\def\no{\noindent}

\def\bin#1#2{\mbox{\scriptsize $\big\{\!\!\!\begin{array}{c} #1 \\ #2
    \end{array}\!\!\!\big\}$\normalsize}}
\renewcommand{\geq}{\geqslant}

\def\sh{{\rm sh}}
\def\hh{{\rm h}}
\def\zmod#1{\mathbb{Z}/#1\mathbb{Z}}
\def\et{{\rm \acute{e}t}}

\def\red{{\rm red}}
\def\alg{{\rm a}}
\def\rig{{\rm rig}}

\def\FA{{\rm (FA)}}
\def\I{\rm I}
\def\II{\rm II}
\def\i{{\rm (i)}}

\def\into{\hookrightarrow}

\def\To{\Rightarrow}
\renewcommand{\tilde}{\widetilde}
\renewcommand{\hat}{\widehat}
\renewcommand{\bar}{\overline}

\newtheorem{counter}[subsubsection]{$\!\!$}
\newtheorem{subcounter}[subsection]{$\!\!$}
\newenvironment{defi}{\begin{counter} \rm {\bf Definition.}}{\end{counter}}
\newenvironment{defis}{\begin{counter} \rm {\bf Definitions.}}{\end{counter}}

\newenvironment{prop}{\begin{counter} {\bf Proposition.}}{\end{counter}}
\newenvironment{lemm}{\begin{counter} {\bf Lemma.}}{\end{counter}}

\newenvironment{theo}{\begin{counter} {\bf Theorem.}}{\end{counter}}
\newenvironment{quot}{\begin{counter} {\bf Theorem \ }}{\end{counter}}

\newenvironment{rema}{\begin{counter} \rm {\bf Remark.}}{\end{counter}}

\newenvironment{exams}{\begin{counter} \rm {\bf Examples.}}{\end{counter}}
\newenvironment{coun}{\begin{counter} \rm {\bf Counter-examples.}}{\end{counter}}
\newenvironment{noth}{\begin{counter} \rm}{\end{counter}}
\newenvironment{proo}{{\flushleft \bf Proof :}}{\hfill $\square$ \vspace{5mm}}

\DeclareMathOperator{\Aut}{Aut}
\DeclareMathOperator{\Hom}{Hom}

\DeclareMathOperator{\Id}{id}

\DeclareMathOperator{\coker}{coker}

\DeclareMathOperator{\Spec}{Spec}
\DeclareMathOperator{\Spf}{Spf}

\DeclareMathOperator{\Rig}{Rig}

\DeclareMathOperator{\PSh}{PSh}

\DeclareMathOperator{\SO}{SO}
\DeclareMathOperator{\ch}{char}

\DeclareMathOperator{\Ass}{Ass}
\DeclareMathOperator{\Ann}{Ann}

\def\cF{{\cal F}}   
\def\cK{{\cal K}}  \def\cM{{\cal M}} \def\cN{{\cal N}}
\def\cO{{\cal O}}

\def\bA{{\mathbb A}}  \def\bG{{\mathbb G}}
  \def\bZ{{\mathbb Z}}

\def\fe{{\mathfrak e}}

\begin{document}

\begin{center}
{\bf \Large Effective models of group schemes}

\bigskip
\bigskip

{\bf Matthieu Romagny}

\medskip

{\em September 28, 2009}
\end{center}

\bigskip

{\def\thefootnote{\relax}
\footnote{ \hspace{-6.8mm}
\\
Matthieu {\sc Romagny}, Institut de Math{\'e}matiques, Th{\'e}orie des Nombres,
Universit{\'e} Pierre et Marie Curie, Case~82,
4, place Jussieu,
F-75252 Paris Cedex 05.

\noindent {\it email}~: romagny@math.jussieu.fr
}}

\begin{center}
\bf{Abstract}
\end{center}

Let $R$ be a discrete valuation ring with fraction field $K$ and $X$ a
flat $R$-scheme. Given a faithful action of a $K$-group scheme $G_K$
over the generic fibre $X_K$, we study models $G$ of $G_K$ acting on
$X$. In various situations, we prove that if such a model $G$ exists,
then there exists another model $G'$ that acts faithfully on $X$. This
model is the schematic closure of $G$ inside the fppf sheaf
$\Aut_R(X)$~; the major difficulty is to prove that it is
representable by a scheme. For example, this holds if $X$ is locally
of finite type, separated, flat and pure and $G$ is finite flat. Pure
schemes (a notion recalled in the text) have many nice properties~:
in particular, we prove that they are the amalgamated
sum of their generic fibre and the family of their finite flat closed
subschemes. We also provide versions of our results in the setting
of formal schemes.


\section{Introduction}

The present paper is interested in the reduction of algebraic
varieties with group action. Let us fix a discrete valuation ring
$R$ with fraction field $K$ and residue field $k$. Algebraic and
arithmetic geometers study all kinds of varieties, or varieties
with additional structures, defined over $K$. In various
situations, these objects have a unique model over $R$ or over
a finite extension~; this is so each time that one has a proper
moduli space for the objects, but not only. Let us mention a few
of these well-known models~: stable models of curves (Deligne,
Mumford), N{\'e}ron models of abelian varieties (N{\'e}ron),
semiabelic pairs as models of principally polarized varieties
(Alexeev), stable maps as models of morphisms from a curve to
a fixed variety (Abramovich, Vistoli). If a group~$G$ acts
faithfully on the $K$-variety and the model satisfies some
unicity property, the action extends to it.

Our concern is, in fact, exclusively in the reduction of the group
action. The point is that even though most of the time the action
of $G$ extends as just indicated, in general the action on the
special fibre is not faithful, and one wishes to consider
other models of~$G$ whose action is better-behaved in reduction.
For typical examples, assume that $R$ has unequal characteristics
$(0,p)$ and $G$ is a finite $p$-group. If $A$ is an abelian
scheme over~$R$, or the N{\'e}ron model of an abelian variety $A_K$
of dimension $g$, such that the $p$-torsion $A_K[p]$ is rational,
then $G=\zmod{p^{2g}}$ acts by translations. This action extends
to~$A$ and, for lack of $p$-torsion points in characteristic~$p$,
the action has a nontrivial kernel on the special fibre.
For another example, consider a smooth pointed curve $(C_K,x_K)$
endowed with a faithful action of $G$ leaving $x_K$ fixed, and
assume that $(C_K,x_K)$ has a stable pointed model $(C,x)$ over~$R$.
We wish to understand the reduction of the action, especially
around the reduction $x_k$. We are led to focus on the orbit
$Z\subset C_k$ of the irreducible component of $x_k$. After
throwing away all components of $C_k$ not in $Z$, we get an
open $R$-curve, and we are asking for the best model for the
induced action of $G$.

In the example above of an abelian scheme $A$, the $R$-group scheme
of $p$-torsion $G'=A[p]$ is the obvious choice of a good model.
We can recover it as follows~: to the action of $G$ is associated
a morphism of $R$-group schemes $G_R\to\Aut_R(A)$, where $G_R$ is
the constant $R$-group scheme defined by $G$. Then $G'$ is the
schematic image of this morphism~; the special properties of
schematic images and closures over a discrete valuation ring
ensure that $G'$ is flat over $R$. In the examples of a N{\'e}ron
model or an open curve, we would like to do the same thing. But
there comes a problem~: these schemes are not proper, and the
automorphism functor is not representable by a scheme or an algebraic
space. Still, it is a sheaf for the fppf topology, and Raynaud has
given a definition of schematic closures in this setting~; but
representability of these closures is by no means obvious, and indeed,
it does not happen in general. The main theorems of this article
prove that these schematic images are often representable by flat
group schemes when we consider actions on {\em pure} schemes, the
notion of purity being a (very) weak version of properness. For
example, faithfully flat $R$-schemes with geometrically irreducible
fibres without embedded components are pure. Here are some
of our most striking results~:

\bigskip

\no {\bf Theorem A.}
{\em {\rm (i)} Let $X$ be an $R$-scheme locally of finite type,
separated, flat and pure. Let $G$ be a proper flat $R$-group scheme
acting on $X$, faithfully on the generic fibre. Let $N$
denote the kernel of the action. Then the schematic image
of $G$ in $\Aut_R(X)$ is representable by a flat group scheme
of finite type $G'$ if and only if $N_k$ is finite. Moreover,
in this case $G'$ is proper.

\smallskip

\no {\rm (ii)} Let $X$ be an affine $R$-scheme, equal to the
spectrum of a ring $A$ such that the map $A\to \Pi\,A/I_\lambda$
to the product of the finite flat quotients of $A$ is universally
injective. Let
$G$ be an $R$-group scheme locally of finite type, flat and pure,
acting on $X$, faithfully on the generic fibre. Then the schematic
image of $G$ in $\Aut_R(X)$ is representable by a flat $R$-group
scheme $G'$. If $G$ is quasi-compact, or affine, or finite, then
$G'$ has the same property.}

\bigskip

When it is representable, we call the schematic image the
{\em effective model} of $G$ for its action on $X$. We have
also versions of these results in the setting of formal schemes.

The affine version in case (ii) is interesting because it applies
not only to rings of finite type, flat and pure (by theorem~B below)
but also for example to rings arising from the completion of smooth
$R$-schemes along a section, and also because the assumptions
made on the group $G$ are very light. Let us now focus on
case~(i). As it turns out, this result does not have much to do
with groups. The crucial facts that govern the proof are the good
properties of $R$-schemes locally of finite type and pure. Such
a scheme $X$ is the amalgamated sum of its generic fibre $X_K$ and
the family of all its closed subschemes finite flat over $R$, the
latter family being schematically dense in a very strong sense.
In fact, we prove~:

\bigskip

\no {\bf Theorem B.}
{\em Assume that $R$ is henselian. Let $X$ be an
$R$-scheme locally of finite type, flat and pure.
Then, the family of all closed subschemes
$Z_\lambda\subset X$ finite flat over $R$ is $R$-universally
schematically dense, and for all separated $R$-schemes $Y$ and
all diagrams in solid arrows
$$
\xymatrix{\amalg Z_{\lambda,K} \ar[r] \ar[d] & X_K \ar[d] \ar[rdd] & \\
\amalg Z_\lambda \ar[r] \ar[rrd] & X \ar@{.>}[rd] & \\
& & Y}
$$
there exists a unique morphism $X\to Y$ making the full
diagram commutative.}

\bigskip

Here also, there is an analogue for formal schemes. Using
Theorem~B, we prove representability results for schematic
images of schemes or formal schemes inside functors of
the type $\Hom_R(X,Y)$. Theorem~A above is essentially an
application of the particular case $X=Y$.

\bigskip

The effective models defined in the present article have been
studied in full detail, for the cyclic group of order $p^2$ in
unequal characteristics, in the recent Ph.D. thesis of
D.~Tossici (\cite{To1}, \cite{To2}). His results provide
more examples of effective models, and show some of their general
features. Related to this work is also the note \cite{Ab} of
Abramovich. There, some group schemes over stable curves
are considered. They are not unrelated with our effective
models, and we plan to compare the two approaches more precisely
in the near future. This will hopefully lead to some new
insights on the reduction of the moduli space of admissible
Galois covers of stable curves (see \cite{BR}). The latter
question is open at the moment, and it was the most important
motivation for the present work.

\subsection{Overview of the article}

Here is a short description of the contents of the article,
together with precise references to the statements of the
main results. In section~\ref{Complements_purity}, we recall
some results on purity and provide some complements.
We prove openness results for some properties of the fibres of
morphisms of finite presentation, flat and pure, that have some
independent interest (theorem~\ref{open_loci_for_pure_morphisms}).
In section~\ref{section_reconstructing} we study schematically
dominant families of morphisms from flat schemes to a fixed
scheme $X$. We prove the density of finite flat closed subschemes
(theorem~\ref{density}) as well as the amalgamme property
(propositions~\ref{cocartesian_diagram_affine_case},
\ref{cocartesian_diagram_general_case},
\ref{cocartesian_diagram_general_formal_case})
which together give the statement of Theorem~B. In the beginning
of section~\ref{Sch_images_inside_Hom_and_Aut} we introduce
schematic images and we prove some useful general results on kernels
for scheme or group scheme actions. Then the stage is set
to prove representability of schematic images in various
situations~: we start with images of schemes inside $\Hom$
functors and then we prove representability of images of groups
in the scheme case
(theorems~\ref{representability_group}
and~\ref{representability_group_proper_case}) and in the formal
scheme case (theorems~\ref{representability_formal_group} and
\ref{representability_formal_group_proper_case}).
Theorem~A is the combination of these results. We also give some
properties enjoyed by the effective model of a finite flat
group scheme (proposition~\ref{various_properties}). Finally,
in section~\ref{Examples} we give some examples. Notably, we
compute explicitly the schematic image in two different cases
of degeneration of torsors under the cyclic group of order $p^2$
in equal characteristic $p>0$ (see~\ref{examples_torsors}).
We observe in particular that for a normal subgroup $H\subset G$,
the effective model of $G/H$ acting on $X/H$ may be different
from $G'/H'$.

\subsection{Notations and conventions}

Everywhere in the paper, we abbreviate the notation of a
discrete valuation ring $R$ with fraction field $K$, residue
field $k$, and chosen uniformizer $\pi$,
by the tuple $(R,K,k,\pi)$. In general, the residue characteristic
is denoted $p\ge 0$. For schemes or morphisms defined over~$R$,
we use subscripts $(\cdot)_K$ and $(\cdot)_k$ to denote the
restrictions to the generic and the special fibre.

When $R$ is complete, we consider also formal $R$-schemes.
A formal scheme $X$ may be identified with a direct system of ordinary
schemes~$X_n$ over the ring $R_n=R/(\pi^n)$. We refer to \cite{BL1}
for basic facts on formal and rigid geometry, and in particular for
the notion of {\em admissible formal blowing-up}. {\em Admissible
formal schemes} in \cite{BL1} are flat formal schemes locally of
finite type. Raynaud's theorem (\cite{BL1}, theorem~4.1)
asserts that there is an equivalence between the category of 
quasi-compact admissible formal $R$-schemes, localised by admissible
formal blowing-ups, and the category of quasi-compact and
quasi-separated rigid $K$-spaces. The $K$-space associated to
a formal scheme $X$ is called its {\em generic fibre} and denoted
$X_\rig$.

\subsection{Acknowledgements}

This paper is the final version of the preprint \cite{Ro},
which it renders outdated and replaces. I wish to thank
Jos{\'e} Bertin and Laurent Moret-Bailly for sharing their
ideas with me when the first outlines of proofs were written.
For various interesting suggestions and remarks,
I am also thankful to Dan Abramovich,
Qing Liu, Michel Raynaud and Dajano Tossici. Finally, thanks are
due to David Eisenbud for pointing out the paper \cite{EH}.

\tableofcontents

\section{Complements on purity} \label{Complements_purity}

%

\subsection{Purity, projectivity and adic topologies}


We first recall some definitions from Gruson-Raynaud \cite{GR}.

\begin{defi}
Let $X\to S$ be a morphism of schemes and $\cM$ be a quasi-coherent
$\cO_X$-module.
\begin{trivlist}
\itemn{i} The {\em relative assassin} of $\cM$ over $S$,
denoted $\Ass(\cM/S)$ is the union over all $s\in S$ of the associated
points $x\in X\otimes k(s)$ of $\cM\otimes k(s)$. If $\cM=\cO_X$ we
set $\Ass(X/S)=\Ass(\cM/S)$.
\itemn{ii} Assume that $X\to S$ is locally of finite type and $\cM$ is
of finite type. We say that $\cM$ is {\em $S$-pure} if the following
condition is satisfied~: for any $s\in S$, if $(\tilde S,\tilde s)$
denotes a henselization of $(S,s)$, then the closure of any point
$\tilde x\in\Ass(\cM\times_S\tilde S/\tilde S)$ meets $X\otimes
k(\tilde s)$. We say that $X$ is {\em $S$-pure} if $\cO_X$ is $S$-pure.
\end{trivlist}
\end{defi}

\begin{exams} \label{examples_pure}
(1) If $X\to S$ is proper, then it is pure.

\no (2) If $X\to S$ is faithfully flat with geometrically
irreducible fibres without embedded components, then it is pure.

\no (3) Let $R$ be a henselian discrete valuation ring and
$X_1=\Spec(R[\varepsilon,x]/(\varepsilon^2,\varepsilon x))$.
Let $X$ be the complement in $X_1$ of the closed point
defined by the ideal $(\pi,\varepsilon,x)$. Then $X$ is not
pure over $R$.
\end{exams}

Here is one of the main results of \cite{GR} (th{\'e}or{\`e}me~3.3.5 in
part~I of {\em loc. cit.})~:

\begin{quot}{\bf (Gruson, Raynaud)} \label{main_thm_from_GR}
Let $A$ be a ring, $B$ an $A$-algebra of finite presentation, $M$ a
$B$-module of finite presentation, flat over $A$. Then $M$ is a
projective $A$-module if and only if it is pure over $A$.
\end{quot}

In what follows, we shall provide some complements on the notion of
purity. In particular, given an $S$-scheme $X$, we will explain the relation
between purity of $X$ and the property that $X$ may have an open covering by
affine schemes with function rings separated for some adic topologies
coming from $S$ (in particular, when $S$ is a local scheme, the maximal-adic
topology). We also give some applications.

\begin{lemm} \label{fppf_pure}
Let $S$ be a scheme and $X,Y$ be $S$-schemes locally of finite type. Let
$f\colon X\to Y$ be an fppf morphism over $S$. Then $Y$ is $S$-pure
if $X$ is $S$-pure. If furthermore $f$ is pure, then the converse holds.
\end{lemm}

\begin{proo}
We may assume that $S$ is a local henselian scheme and since the
locus of the base where a map is pure is open (\cite{GR},~I.3.3.8),
it is enough to test purity at the closed point $s\in S$. Now
let $y\in\Ass(Y/S)$. Choose some associated point $x\in X_y$ so
$x\in\Ass(X/S)$. Then there exists $a\in X_s$ meeting the closure
of $x$, so $f(a)$ meets the closure of $y$. So $Y$ is $S$-pure.

Conversely, assume that $f$ is pure and let $x\in\Ass(X/S)$ and
$y=f(x)$. Thus $x\in\Ass(X/Y)$ and $y\in\Ass(Y/S)$
(\cite{GR},~I.3.2.4). Since
$Y$ is $S$-pure, the closure of $y$ meets $Y_s$ at some point $b$. Let
$(\tilde Y,\tilde b)$ be a henselization of $(Y,b)$, let $\tilde X=X\times_Y
\tilde Y$, and $\tilde x=(x,\tilde b)\in \tilde X$ so that $\tilde
x\in\Ass(\tilde X/\tilde Y)$ by \cite{GR},~I.3.2.3. Thus the closure of
$\tilde x$ inside $\tilde X$ meets $\tilde X_{\tilde b}$ at a point $\tilde
a$. The image of $\tilde a$ in $X$ lies in the closure of $x$ and above
$\tilde b$, thus in $X_s$. Therefore $X$ is $S$-pure.
\end{proo}

%
%

\begin{defi}
Let $n\ge 1$ be an integer.
We say that a morphism of schemes $X\to S$ is {\em of type \FA$_n$} is every
set of $n$ points of $X$ whose images lie in an affine open set of $S$ lie
in an affine open set of $X$. We say that $X\to S$ is {\em of type \FA}
if it is of type \FA$_n$ for all $n\ge 1$.
\end{defi}

\begin{lemm} \label{FA_and_limit}
Assume that $S$ is affine. Let $X\to S$ be of finite presentation and
of type \FA$_n$. Then there exists a scheme $S_0$ of finite type over
$\bZ$ and a morphism $S\to S_0$, an $S_0$-scheme $X_0$ of finite
presentation of type \FA$_n$, such that $X\simeq X_0\times_{S_0}S$.
\end{lemm}

\begin{proo}
Since $S$ is affine and $X\to S$ is quasicompact, to say that
$X\to S$ is of type \FA$_n$ means that there exists a finite cover by open
affine schemes $U_i$ ($1\le i\le m$) such that $\amalg (U_i)^n\to X^n$ is
surjective. Thus using \cite{EGA}~IV.8.10.5~(vi), there exists a scheme
$S_0$ of finite type over $\bZ$, an $S_0$-scheme $X_0$ of finite presentation,
and an open cover $U_{0,i}$ of $X_0$ such that $U_i\simeq
U_{0,i}\times_{S_0}S$ for all $i$, $X\simeq X_0\times_{S_0}S$, and $\amalg
(U_{i,0})^n\to (X_0)^n$ is surjective.
\end{proo}

In the next lemma, we relate the notion of purity for a scheme
over a noetherian henselian local ring $(R,m)$ with the property of
separation of the function rings with respect to the $m$-adic topology.
We will say that an $R$-algebra $A$ is {\em strongly separated for the
$m$-adic topology} if and only if for all prime ideals $q\subset m$,
the ring $A/qA$ is separated for the $m/q$-adic topology.

\begin{lemm} \label{purity_versus_adic_separation}
Let $R$ be a noetherian henselian local ring with maximal ideal $m$.
Let $X$ be a scheme locally of finite type and flat over $R$. Consider the
following conditions~:
\begin{trivlist}
\itemn{i} $X$ is $R$-pure.
\itemn{ii} $X$ has an open covering by affine schemes whose function
algebras are free $R$-modules.
\itemn{iii} $X$ has an open covering by affine schemes whose function
algebras are strongly $m$-adically separated.
\end{trivlist}
Then, we have {\rm (ii)} $\To$ {\rm (iii)} $\To$ {\rm (i)}.
Assume moreover that $R$ is a discrete valuation ring and $X$ is
quasicompact and of type \FA$_{n+1}$, where $n$ is the number of
associated points of the generic fibre. Then all three conditions
{\rm (i)}, {\rm (ii)}, {\rm (iii)} are equivalent. Furthermore, we
may choose an open covering $\{U_i\}$ as in {\rm (ii)}-{\rm (iii)}
so that all intersections $U_i\cap U_j$ are $R$-pure again.
Finally, the $R$-module $H^0(X,\cO_X)$ is free.
\end{lemm}

\begin{proo}
The fact that (ii) implies (iii) is clear since any free $R$-module
is strongly separated for the $m$-adic topology. Let us check
that (iii) implies (i). Let $x\in \Ass(X/R)$ and $U=\Spec(A)$ be an
open affine containing $x$, with $A$ strongly $m$-adically separated.
Let $p\subset A$ (resp. $q\subset R$) be the prime ideal
corresponding to $x$ (resp. the image of $x$ in $S$) and let
$k(q)=A_q/qA_q$ be the residue field of $q$. If the closure of $x$
in $U$ does not meet the special fibre, there exist $u\in p$ and
$v\in m$ such that $1=u+v$. But by assumption, there is $a\in A$
such that the image of $p$ in $A\otimes k(q)$ is the annihilator
$\Ann_{A\otimes k(q)}(a)$. Hence there exists $s\in R\setminus q$
such that $sua\in qA$. In the ring $A/qA$, we get $sa=sav=sav^n$
for all $n\ge 1$ hence $sa$ lies in $\cap_{n\ge 0}\,(m/q)^n(A/q)$.
The latter intersection is zero by assumption, hence $sa=0$ in
$A/q$ and $a=0$ in $A\otimes k(q)$. This is impossible. By
contrapositive, $X$ is pure.

We now prove that under the additional assumptions,
we have (i)~$\To$~(ii). Call $x_1,\dots,x_n$ the associated points
of the generic fibre of $X$.
By purity, for each $i$ the closure of $x_i$ meets the closed fibre in
at least one point $x'_i$. Since it is assumed that $X\to \Spec(R)$ is
of type \FA$_{n+1}$, for each $x\in X$ we may find an open affine
$U_x=\Spec(A)$ containing $x,x'_1,\dots,x'_n$. Obviously $U_x$ is
$R$-pure, so it follows from~\ref{main_thm_from_GR} that $A$ is a projective
$R$-module, i.e. a free $R$-module since $R$ is a principal ideal domain.
Since $X$ is quasicompact, we can extract from $\{U_x\}$ a finite
open cover, and since each of them contains $x'_1,\dots,x'_n$,
the intersections $U_i\cap U_j$ are $R$-pure.

Finally we prove that $H^0(X,\cO_X)$ is free.
Let $U_i=\Spec(A_i)$ be an open covering by affine
schemes whose function algebras are free $R$-modules. Since $X$ is
quasi-compact, finitely many of the $U_i$ are sufficient. Since a
submodule of a free $R$-module is free, the injection
$H^0(X,\cO_X)\hookrightarrow \Pi\,H^0(U_i,\cO_{U_i})$ gives the
desired result.
\end{proo}

\begin{rema} \label{rema_pure_affine_case}
The special case where $X$ is affine of finite type and flat over
a discrete valuation ring will be useful later in the paper. In this
case, the proof above shows that $X$ is pure if and
only if $\Gamma(X,\cO_X)$ is a free $R$-module, if and only
if $\Gamma(X,\cO_X)$ is separated for the $\pi$-adic topology.
\end{rema}


We now point out some features of pure schemes over a discrete valuation
ring, and in particular a relation between purity and the topology of
the neighbourhoods of the special fibre. Note that the notions of
schematic density and schematic dominance will receive a more complete
treatment in section~\ref{section_reconstructing}~; we refer to it
for more details.

\begin{lemm} \label{neighbourhood_sch_dense}
Let $(R,K,k,\pi)$ be a discrete valuation ring.
The following properties hold.
\begin{trivlist}
\itemn{1} Let $f:Z\to X$ be a morphism of $R$-schemes with $X$
flat over $R$. Then $f$ is schematically dominant if and
only if $f_K$ is schematically dominant.
\itemn{2} Let $X$ be an $R$-scheme locally of finite type and pure.
Then any open neighbourhood of the closed fibre $X_k$ is
schematically dense in $X$. If moreover $X$ is flat over $R$,
then such a neighbourhood is $R$-universally schematically dense.
\itemn{3} Let $X,Y$ be $R$-schemes of finite type with $X$ pure and
$Y$ separated. Let $\hat R$, $\hat X$, $\hat Y$ be the $\pi$-adic
formal completions of $R$, $X$, $Y$. Then, the completion map
$$
\Hom_R(X,Y) \to \Hom_{\hat R}(\hat X,\hat Y) \quad,\quad f\mapsto \hat f
$$
is injective.
\end{trivlist}
\end{lemm}


\begin{proo}
(1) This is clear, since $X_K$ is schematically dense in $X$.

\no (2)
Let $U$ be an open neighbourhood of the closed
fibre $X_k$. To prove that $U$ is schematically dense, we may
replace~$R$ by its henselization and hence assume that $R$ is
henselian. Then it is enough to prove that $\Ass(X)\subset U$.
If $x\in \Ass(X)$, then by \cite{EGA}~IV.3.3.1 it is an associated
point in its fibre $X_s$, where $s$ is the image of $x$ in $s$.
Since $X$ is pure, the closure of $x$ meets $X_k$, hence it
meets $U$, so $x\in U$ and we are done. If moreover $X$ is flat
over $R$, then using point~(1) we see that $U_K$ is schematically
dense in $X_K$. Since $U_k=X_k$ is schematically dense in $X_k$
and~$U$ is flat over $R$, it is $R$-universally schematically dense.

\no (3) Let $f,g:X\to Y$ be such that $\hat f=\hat g$. In order to
prove that $f=g$, we may pass to the henselization of $R$ and hence
assume that $R$ is henselian. By \cite{EGA}~I.10.9.4, there is
an open neighbourhood $U\subset X$ of $X_k$ where $f$ and $g$ are
equal. It follows from (1) that $U$ is schematically dense in $X$.
Since $f=g$ on a schematically dense open subscheme of
$X$, we get $f=g$ on $X$.
\end{proo}

\begin{rema}
Point (2) of this lemma allows to compare pure schemes with other
schemes by looking at ''how close'' a scheme is to its special
fibre. If we arrange $R$-schemes by increasing distance to their
special fibre, we have $k$-schemes, then formal $R$-schemes,
then pure $R$-schemes, then general $R$-schemes.
\end{rema}

\begin{lemm} \label{msep}
Let $X\to S$ be a morphism of schemes. Assume that $X$ is locally noetherian
and $S$ is affine. Let $s\in S$ and let $p\subset \Gamma(S,\cO_S)$ be the
corresponding ideal. Then there is an open neighbourhood of the fibre $X_s$
that is covered by affine schemes whose function ring is separated for the
$p$-adic topology.
\end{lemm}

\begin{proo}
Let $x\in X_s$ and let $U_1=\Spec(A_1)$ be an affine neighbourhood
with $A_1$ noetherian. Let $m\subset A_1$ be the prime ideal
corresponding to $x$, so $pA_1\subset m$. Let
$I_1=\cap_{n\ge 0}\,p^nA_1$. Since $\cO_{X,x}$ is local noetherian,
it is separated for the $p$-adic topology, hence $I_1$ lies in the
kernel of the localization morphism
$A_1\to\cO_{X,x}$. Since $I_1$ is finitely generated, there is $s_1\in
A_1\setminus m$ such that $s_1I_1=0$. In other words, if we set
$A_2=A_1[1/s_1]$ then $I_1$ maps to $0$ under $A_1\to A_2$. By
induction, after we have defined $A_r$, we let
$I_r=\cap_{n\ge 0}\,p^nA_r$, we argue that there is
$s_r\in A_r\setminus m$ such that $s_rI_r=0$, and we define
$A_{r+1}=A_r[1/s_r]$. Because $A_1$ is noetherian, the increasing
sequence of ideals $K_r=\ker(A_1\to A_r)$ must stabilise at
some~$\rho$. One checks that $I_\rho=0$, that is, $A_\rho$ is separated
for the $p$-adic topology.
\end{proo}

\subsection{Application to the fibres of morphisms}

We now mention an application of these results to the study of the fibres of
morphisms of schemes. Namely, one can weaken the assumptions in some theorems
of \cite{EGA}~IV,~\S~12.2, by requiring purity instead of properness.

\begin{theo} \label{open_loci_for_pure_morphisms}
Let $f\colon X\to S$ be of finite presentation, flat and pure, and
let $n\ge 1$ be an integer. Then the following sets are open in $S$~:
\begin{trivlist}
\itemn{i} The locus of points $s\in S$ such that the fibre $X_s$ is
geometrically reduced.
\itemn{ii} The locus of points $s\in S$ such that the geometric fibre
$X_{\bar s}$ is reduced with less than $n$ connected components.
\itemn{iii} The locus of points $s\in S$ such that the geometric fibre
$X_{\bar s}$ is reduced and has less than $n$ irreducible components.
\end{trivlist}
\end{theo}

\begin{proo}
The assertions to be proven are local on $S$ so we may assume
$S=\Spec(R)$ affine.
By limit arguments using
\cite{GR}, corollaire 3.3.10, and other usual results of
\cite{EGA}~IV,~\S\S~8--11, we reduce to the case where~$R$
is noetherian.

Let $P$ be one of the properties {\em reduced}, {\em reduced
and connected}, or {\em integral}.
The loci we are interested in are constuctible so it is enough
to prove that they are stable under generization. By
\cite{EGA}~II.7.1.7 one reduces to $R=(R,K,k,\pi)$ equal to a
discrete valuation ring, which we may assume henselian.
Then we assume
that the closed fibre has the geometric property $P$, and we have
to prove that the generic fibre has it also. For this it is
enough to prove that for all finite field extensions $L/K$, the
scheme $X\otimes_K L$ has property~$P$. Replacing $R$ by its
integral closure in $L$ we reduce to $K=L$. We now consider the three
cases separately.

\medskip

\no (i) By lemma~\ref{msep}, there is an open neighbourhood $U$ of the
special fibre of $X$ that is covered by open affine subschemes with
function ring separated for the $\pi$-adic topology, i.e. pure.
By lemma~\ref{neighbourhood_sch_dense} this $U$ is universally
schematically dense so if~$U$ has reduced generic fibre then $X$
also. Therefore we may replace $X$ by $U$ and hence assume that~$X$ is
covered by pure open affine subschemes. Let $V=\Spec(A)$ be such an
open affine, it is enough to prove that $A$ is reduced. Since $A$ is
separated for the $\pi$-adic topology and has no $\pi$-torsion, if $x$
is a nonzero nilpotent we may assume that $x\not\in \pi A$. But then
we have a contradiction with the fact that $A_k$ is reduced. So $A$ is
reduced.

\medskip

\no (ii) From (i) we know that $X_K$ is reduced. Then we may as
in (i) reduce to the case where $X$ is covered by pure open affine
subschemes. We shall prove that the number of connected components of
$X_K$ is less than that of $X_k$.
Let $B=H^0(X,\cO_X)$. From the injection
$B_k\hookrightarrow H^0(X_k,\cO_{X_k})$ we learn that $B_k$ is
reduced. This, together with an easy calculation, proves that
the idempotents of~$B$ and those of $B_K$ are the same. So $X_K$
and $X$ have the same number of connected components~; call it
$u$. Then $B$ splits as a product of rings
$B_1\times\dots\times B_u$, with $B_i\ne 0$ for $i=1,\dots,u$.
Since $B$ is a free $R$-module
(lemma~\ref{purity_versus_adic_separation}), each of the $B_i$
is free, and hence $B_{i,k}\ne 0$. Hence $B_k$ has at least $2^u$
idempotents, so $X_k$ has at least $u$ connected components.

\medskip

\no (iii) From (ii) we know that $X_K$ is reduced and connected.
It is enough to prove that for any irreducible component $W$ of
$X_k$, there is a unique irreducible component $Z$ of $X_K$ whose
closure in $X$ contains $W$. For this, we may remove from $X$ all
irreducible components $W'\ne W$ of $X_k$ and all irreducible
components $Y$ of $X$ that do not contain $W$ (they are closures in
$X$ of irreducible components of $X_K$). Hence, we may assume that
$X_k$ is integral, and we have to prove that $X_K$ is integral also.
We may as in (i) reduce to the case where $X$ is covered by pure open
affine subschemes. It is then enough to prove that all such open
affines $V=\Spec(A)$ are integral. But if $xy=0$ in $A$, and $x,y$ are
nonzero, we may as in (i) assume that they do not belong to $\pi
A$. Then this contradicts the fact that $A_k$ is integral.
\end{proo}

\begin{coun}
Obviously the corollary does not extend to all properties listed in
\cite{EGA}, IV, \S~12.2. We give counter-examples for some of them.
Let $(R,K,k,\pi)$ be a discrete valuation ring.
\begin{trivlist}
\itemn{1} {\em geometrically connected}. Let $A=R[t]/(t^2-\pi t)$ and
$X=\Spec(A)$. Then $X_k$ is geometrically connected but $X_K$ has two
connected components.
\itemn{2} {\em geometrically pointwise integral}.
Let $A=R[e,x,y]/I$ where $I$ is the ideal generated by the four
elements $xy$, $e^2-e+\pi$, $(1-e)x-\pi x$, $ey-\pi y$.
Let $X=\Spec(A)$. Then $X_k$ is geometrically
pointwise integral (with two connected components), but $X_K$ is not, for it
is geometrically connected and $A_K$ has zerodivisors $x,y$.
\itemn{3} {\em smooth, geometrically normal, etc}. Let $X$ be a flat finite
type $R$-scheme with geometrically integral fibres without embedded
components. Let $U$ be the complement in $X$ of the singular locus of
$X_k$. Then $U$ is again pure over $R$, with smooth special fibre, but the
generic fibre can be chosen to haved arbitrary singularities.
\end{trivlist}
\end{coun}

\section{Reconstructing a scheme from flat closed subschemes}
\label{section_reconstructing}

In this section, we consider two types of situations~:
\begin{trivlist}
\itemn{\I} {\em Ordinary}~: a discrete valuation ring $(R,K,k,\pi)$ and an
$R$-scheme $X$ with a family of morphisms of $R$-schemes
$Z_\lambda\to X$ indexed by a set $L$.
\itemn{\II} {\em Formal}~: a complete discrete valuation ring $(R,K,k,\pi)$
and a formal $R$-scheme $X$ with a family of morphisms of formal
$R$-schemes $Z_\lambda\to X$ indexed by a set $L$.
\end{trivlist}

\medskip

Most of the time, we write this family as a single morphism
$f:\amalg\, Z_\lambda\to X$. After some generalities in subsection
\ref{sdm}, we specialise in subsection~\ref{gaffs} to the case
where $f$ is the family of all
(formal) closed subschemes of $X$ finite flat over $R$. The
general theme is to find some conditions under which $X$ is the
amalgamated sum of its generic fibre $X_K$ and the subschemes
$Z_\lambda$ along the subschemes $Z_{\lambda,K}$ (in the formal
case, the generic fibres are the rigid analytic spaces $X_\rig$
and $Z_{\lambda,\rig}$). As a matter of notation, when no confusion
seems possible, we will allow ourselves a slight abuse by
maintaining the letter $f$ to denote the restriction
$Z_{\lambda_0}\hookrightarrow \amalg\, Z_\lambda\to X$, for a given
$\lambda_0\in L$. For example, we will write $f_*\cO_{Z_\lambda}$
instead of $(f_{|Z_\lambda})_*\cO_{Z_\lambda}$.

\subsection{Schematically dominant morphisms} \label{sdm}

We will need various notions of dominant
morphisms~; see also \cite{EGA}~IV.11.10.

\begin{defis} \label{defis_dominance}
Let $f:\amalg\, Z_\lambda\to X$, $\lambda\in L$, be a family of
morphisms of $R$-schemes.
\begin{trivlist}
\itemn{1} If $X$ is affine, $f$ is called {\em affinely dominant}
if the intersection of the kernels of the maps $\Gamma(X,\cO_X)\to
\Gamma(Z_\lambda,\cO_{Z_\lambda})$ is $0$. If $X$ is arbitrary,
$f$ is called {\em weakly schematically dominant} if there
exists a covering of $X$ by open affine subschemes $U_i$ such that
$f^{-1}(U_i)\to U_i$ is affinely dominant for all $i$.
\itemn{2} The map $f$ is called {\em schematically dominant} if the
intersection of the kernels of the maps of sheaves $\cO_X\to
(f_\lambda)_*\cO_{Z_\lambda}$ is $0$, or equivalently, if for {\em
  all} open affine subschemes $U\subset X$, the map $f^{-1}(U)\to U$
is affinely dominant.
\end{trivlist}
If one of these properties is true after any base change $R\to R'$, we
say that it is true universally.
\end{defis}

For example, the family of maps
$\amalg_{n\ge 0}\,\Spec(R/\pi^n)\to \Spec(R)$ is affinely dominant,
hence weakly schematically dominant, but not schematically dominant.

If $X$ is affine, it is equivalent to say that $f$ is affinely
dominant or that for any two morphisms $u,v\colon U\to X'$
to an affine $R$-scheme $X'$, $f\circ u=f\circ v$ implies
$u=v$. If $X$ is arbitrary, it is equivalent to say that $f$ is
schematically dominant or that for
any open set $U\subset X$, and any two morphisms $u,v\colon U\to X'$
to a separated $R$-scheme $X'$, if the compositions of $u$ and $v$
with the restriction $f^{-1}(U)\to U$ are equal, then $u=v$.
In the case where each
$f_{|Z_\lambda}$ is an immersion, this gives the notion of a
{\em schematically dense} family of subschemes.

If we consider a family of morphisms of formal $R$-schemes
$f:\amalg\, Z_\lambda\to X$, $\lambda\in L$, the same definitions
and remarks apply word for word.

In the sequel, we will meet one particular case where weakly
schematically dominant morphisms are schematically dominant. In order
to explain this, we recall
the following standard notation : if $I,J$ are ideals in a ring $A$,
we write $(I:J)_A$ or simply $(I:J)$ for the ideal of elements $a\in A$
such that $aJ\subset I$, and we write $(I:J^\infty)$ for the
increasing union of the ideals $(I:J^n)$. The following definition
applies in the case of schemes or formal schemes.



\begin{defi} \label{torsion_bounded}
We say that {\em the torsion in $f_*\cO_{Z_\lambda}$ is
bounded uniformly in $\lambda$} if and only if for all $U\subset X$
open, for all $t\in \cO_X(U)$, there exists an integer $c\ge 1$
such that for all $\lambda\in L$, we have $(0:t^\infty)=(0:t^c)$
as ideals of $(f_*\cO_{Z_\lambda})(U)$.
\end{defi}

\begin{lemm} \label{criterion_bounded_torsion}
Let $f:\amalg\, Z_\lambda\to X$, $\lambda\in L$, be a family of
morphisms of $R$-schemes or formal $R$-schemes. Assume that
either $L$ is finite, or the torsion in $f_*\cO_{Z_\lambda}$ is
bounded uniformly in $\lambda$. Then $f$ is schematically dominant if
and only if it is weakly schematically dominant.
\end{lemm}

\begin{proo}
Only the {\em if} part needs a proof. Let $U=\Spec(A)$ in the scheme
case, resp. $U=\Spf(A)$ in the formal scheme case, be an open affine
such that $f^{-1}(U)\to U$ is affinely dominant. Let
$B_\lambda=(f_*\cO_{Z_\lambda})(U)$,
$\varphi_\lambda:A\to B_\lambda$ the map corresponding to $f_\lambda$,
and $I_\lambda=\ker(\varphi_\lambda)$. The intersection of the ideals
$I_\lambda$ is zero and we have to prove that for all $t\in A$, the
intersection of the kernels of the maps 
$\varphi_\lambda[1/t]: A[1/t]\to B_\lambda[1/t]$
is zero. Let $a$ be in this intersection. Clearly it is enough to
take $a\in A$. For all $\lambda$ there is an integer
$c_\lambda\ge 0$ such that $t^{c_\lambda}\varphi_\lambda(a)=0$. If the
torsion in $f_*\cO_{Z_\lambda}$ is bounded uniformly in $\lambda$,
there is an integer $c$ such that for all $\lambda$ we have
$t^c\varphi_\lambda(a)=0$. If $L$ is finite, this is also true with
$c=\sup\{c_\lambda\,,\,\lambda\in L\}$. It follows that $t^ca$ is in the
intersection of the $I_\lambda$, hence zero by assumption. Thus $a=0$
in $A[1/t]$.
\end{proo}

We now use more specifically the properties of flat modules over the
discrete valuation ring~$R$.
The first lemma below is stated as a useful
observation to keep in mind. Then we continue with some properties
of schemes dominated by flat families.

\begin{lemm} \label{equivalent_conditions_univ_injective}
For a morphism of $R$-modules $u:M\to N$ with $N$ flat, the following
conditions are equivalent~:
\begin{trivlist}
\itemn{1} $u$ is universally injective.
\itemn{2} $u$ is injective and $u_k$ is injective.
\itemn{3} $u$ is injective and $\coker(u)$ is flat.
\end{trivlist}
If $N$ is a direct product of flat modules $N_\lambda$, $\lambda\in L$, and
we denote by $I_\lambda$ the kernel of $M\to N_\lambda$, these
conditions are also equivalent to~:
\begin{trivlist}
\itemn{4} $\underset{\lambda\in L}{\cap}\, I_\lambda=0$ and
$\underset{\lambda\in L}{\cap}\,I_{\lambda,k}=0$.
\end{trivlist}
\end{lemm}

\begin{proo}
This is classical.
\end{proo}

The main point of the following result is to say that $X$ satisfies the
property of the amalgamated sum of $X_K$ and the $Z_\lambda$ along their
respective generic fibres, for morphisms to affine $R$-schemes $Y$.

\begin{prop} \label{cocartesian_diagram_affine_case}
Let $f:\amalg\, Z_\lambda\to X$ be a family of morphisms of $R$-schemes
with $Z_\lambda$ flat over~$R$, for all $\lambda\in L$. Assume moreover
that we are in one of the following cases.
\begin{trivlist}
\itemn{i} $X$ has a covering by open affine schemes $U_i$ whose function
algebras are $\pi$-adically separated and the restriction of $f_k$
to $f^{-1}(U_i)_k$ is affinely dominant.
\itemn{ii} $X$ is locally noetherian and $f_k$ is schematically dominant.
\end{trivlist}
Then the following properties hold.
\begin{trivlist}
\itemn{1} $X$ is flat over $R$.
\itemn{2} $f$, equivalently $f_K$, is weakly schematically dominant
(in case (ii) one needs to assume also that $X$ is locally of finite
type and pure).
\itemn{3} For all affine $R$-schemes $Y$ and all
diagrams in solid arrows
$$
\xymatrix{\amalg Z_{\lambda,K} \ar[r] \ar[d] & X_K \ar[d] \ar[rdd] & \\
\amalg Z_\lambda \ar[r] \ar[rrd] & X \ar@{.>}[rd] & \\
& & Y}
$$
there exists a unique morphism $X\to Y$ making the full diagram
commutative.
\end{trivlist}
\end{prop}

\begin{proo}
Observe that after we have proven that $X$ is flat, in order to
prove the amalgamated sum property to affine schemes, since $X$
is flat and $Y$ is separated, the map $g\colon X\to Y$ is unique
if it exists. Thus we may define it locally on $X$ and glue.
It follows that all assertions to be established are local.

In case (i) we are immediately reduced to the situation where
$X=\Spec(A)$ with $A$ separated for the $\pi$-adic
topology. We keep the notations of the proof of
lemma~\ref{criterion_bounded_torsion} and we also set $B=\Pi\,B_\lambda$,
$\varphi=\Pi\,\varphi_\lambda$ and $I=\ker(\varphi)$. From the injection
$A/I\hookrightarrow B$ it follows that $A/I$ has no $\pi$-torsion
hence is flat over $R$. If $a\in I$, then since $\varphi_k$ is
injective, there exists $a_1\in A$ such that $a=\pi a_1$. Since
$A/I$ has no $\pi$-torsion, $a_1$ himself lies in $I$, and by
induction we obtain $a\in\cap\pi^nA$. So $a=0$ by the assumption
on $A$. This proves that $A$ is torsion-free, hence flat over~$R$,
and also that $f$ is weakly schematically dominant. Now we have
a diagram with all morphisms injective~:
$$
\xymatrix@R=16pt{B \ar[r] & B_K \\
A \ar[r] \ar[u] & \ar[u] A_K}
$$
Obviously, in order to prove the amalgamated sum property for
maps to affine schemes, it is enough to show that $A$ is isomorphic
to the fibred product $A_K\times_{B_K} B$. Since $A$ is separated
for the $\pi$-adic topology, a nonzero element in
$B\cap A_K$ may be written $a/\pi^d$ with $a\in A$ and $d\in\bZ$
minimal, such that there exists $b\in B$ with $a=\pi^d b$ in $B$.
If $d\geq 1$, reducing modulo $\pi$ we find that the image of
$a$ vanishes in $B_k$. Since $A_k\to B_k$ is injective, it follows
that $a\in \pi A$, and this contradicts the minimality of~$d$.
Hence $d\le 0$, so $a/\pi^d\in A$ and we are done.

In case (ii), in order to prove flatness it is enough to look at
points of the special fibre $X_k$. By lemma~\ref{msep}, such a point
has an affine neighbourhoood $\Spec(A)$ with $A$ separated for the
$\pi$-adic topology. From case (i) follows that $X$ is flat.
Also, in this way we have found a neighbourhood $U$ of the special
fibre which is covered by open affine schemes whose function
algebras are $\pi$-adically separated. From case (i) follows that
the restriction of $f$ to $f^{-1}(U)$ is weakly schematically dominant.
So if $X$ is locally of finite type and pure, $U$ is schematically
dominant in $X$ by lemma~\ref{neighbourhood_sch_dense}, hence
$f$ itself is weakly schematically dominant.
Finally, to prove the amalgamated sum property, it is enough to define
$g$ in a neighbourhood of all closed points $x\in X_k$. By
lemma~\ref{msep} we may choose a neighbourhood $\Spec(A)$ where
$A$ is $\pi$-adically separated. Then we are reduced to case (i).
\end{proo}

It is possible to formulate an analogue of the amalgamated sum
property for formal schemes finite type, using the definition of the
generic fibre as a rigid analytic $K$-space as in~\cite{BL1}. Since
we have to impose the assumption of finite type, the direct formal
analogue of the affine version~\ref{cocartesian_diagram_affine_case}
is not relevant. Hence we will content ourselves with a statement of
the properties needed in order to
prove~\ref{cocartesian_diagram_general_formal_case}.

\begin{prop} \label{cocartesian_diagram_affine_formal_case}
Assume that $R$ is complete. Let $f:\amalg\, Z_\lambda\to X$
be a family of morphisms of formal $R$-schemes locally of
finite type, with $Z_\lambda$ flat over $R$ for all $\lambda\in L$,
such that $f_k$ is schematically dominant. Then,
\begin{trivlist}
\itemn{1} $X$ is flat over $R$.
\itemn{2} $f$ (equivalently $f_K$) is weakly schematically dominant.
\end{trivlist}
\end{prop}

\begin{proo}
(1) We may restrict to an open affine formal
subscheme $\Spf(A)$. Then $A$ is $\pi$-adically separated and the
arguments of the proof of point (1) in
proposition~\ref{cocartesian_diagram_affine_case} carry on.

\no (2) The arguments are the same as in point (2) in
proposition~\ref{cocartesian_diagram_affine_case}.
\end{proo}

In the sequel of the paper, we will be mainly interested in the case where $L$
is infinite.
Concerning the case where $L$ is finite (this is essentially the case where
$L$ has just one element, for, one may consider $Z=\amalg\, Z_\lambda$),
the following property is still worth recording~:

\begin{prop}
Let $S$ be a scheme and let $f\colon Z\to X$ be a morphism
of flat $S$-schemes of finite presentation. Assume that $X$ is pure.
Let $S_0\subset S$ be the locus of points $s\in S$ such that $f_s$
is schematically dominant, $X_0=X\times_S S_0$, $Z_0=Z\times_S S_0$.
Then $S_0$ is open in $S$ an $f_{|Z_0}\colon Z_0\to X_0$
is $S_0$-universally schematically dominant.
\end{prop}

\begin{proo}
As in the proof of theorem~\ref{open_loci_for_pure_morphisms},
one reduces to the case where $S$ is the spectrum of a henselian
discrete valuation ring $R$ with uniformizer $\pi$, and $f_k$
is schematically dominant.
By lemma~\ref{msep}, there is an open neighbourhood $U$ of the
special fibre of $X$ that is covered by open affine schemes whose
function algebras are $\pi$-adically separated. By
\ref{cocartesian_diagram_affine_case}(2) and
\ref{criterion_bounded_torsion}, the restriction of $f$ to $U$ is
schematically dominant. Since $U$ is schematically dense in $X$ by
lemma~\ref{neighbourhood_sch_dense}, then $f$ is
schematically dominant.
The fact that $f_{|Z_0}\colon Z_0\to X_0$ is $S_0$-universally
schematically dominant is a consequence of \cite{EGA}~IV,~11.10.9.
\end{proo}

\subsection{Glueing along the finite flat subschemes}
\label{gaffs}

We continue with the ordinary (\I) and formal (\II) situations
presented at the beginning of section~\ref{section_reconstructing}.
From now on, the family $Z_\lambda$ will always be the
family of all closed subschemes of $X$ in case (\I), resp. closed
formal subschemes of $X$ in case (\II), that are finite and flat
over $R$. We denote this family by $\cF(X)$.
Under some mild conditions, we will
prove that this family is $R$-universally schematically dense in
$X$ and we will improve proposition~\ref{cocartesian_diagram_affine_case}
by extending the amalgamated sum property to morphisms to
arbitrary separated (formal) schemes $Y$.

We keep the notation $f\colon \amalg Z_\lambda\to X$ for the canonical
morphism induced by the inclusions $Z_\lambda\subset X$.
Note that $\cF(X)$ is naturally an inductive system,
if we consider it together with the closed immersions
$Z_\lambda\hookrightarrow Z_\mu$. Moreover, we can define the union
of two finite flat closed subschemes by the intersection of the
defining ideals~; this is again a finite flat closed
subscheme. In this way, we see that $\cF(X)$ is filtering.

Let us start our programme. We start with a well-known property.

\begin{lemm} \label{closed_points_attained}
Consider one of the two situations~:
\begin{trivlist}
\itemn{\I} $X$ is an $R$-scheme locally of finite type. Assume that
$X$ is flat over $R$, or more generally that $X_\red$ is flat over $R$.
\itemn{\II} $R$ is complete and $X$ is a flat formal $R$-scheme
locally of finite type.
\end{trivlist}
Then $\cF(X)$ attains all the closed points of $X_k$. In case
{\rm (\I)} the converse is true~: if $\cF(X)$ attains all the
closed points of $X_k$ then $X_\red$ is flat over $R$.
\end{lemm}

\begin{proo}
In case (\I), first note that $X_\red$ is flat if and only if no
irreducible component of $X$ is included in the special fibre.
Hence if $X_\red$ is flat, for each closed point $x\in X_k$,
there is an irreducible component $W\subset X$ at $x$ that is
not contained in $X_k$. Then the claim follows from proposition
10.1.36 of \cite{Liu} applied to $W$. Conversely if $X_\red$ is
not flat then there is an irreducible component included in the
special fibre, and it is clear that this component contains at
least one point not lying on any $Z\in\cF(X)$. In case~(\II)
this is just \cite{BL1}, proposition 3.5.
\end{proo}

For the sequel, a crucial ingredient is a theorem of Eisenbud
and Hochster (see~\cite{EH}) which we recall for convenience~:

\begin{quot} {\bf (Eisenbud, Hochster)} \label{thm_EH}
Let $A$ be a ring, and let $P$ be a prime ideal of $A$. Let $\cN$ be a
set of maximal ideals $m$ such that $A_m/P_m$ is a regular local ring,
and such that
$$
\bigcap_{m\in \cN}\,m=P \ .
$$
If $M$ is a finitely generated $P$-coprimary module annihilated by $P^e$,
then
$$
\bigcap_{m\in \cN}\,m^eM=0 \ .
$$
\end{quot}

As a preparation for the proof of theorem~\ref{density} below,
we first establish
a lemma. We refer to Bruns-Herzog~\cite{BH} for more details on the
following notions. Let $(A,m)$ be a noetherian local ring of
dimension $r$, and write $\lg_A(M)$ or simply $\lg(M)$ for the
length of an $A$-module~$M$. For an arbitrary ideal of definition
$q\subset A$, one defines the Hilbert-Samuel multiplicity $\fe(q)$
as the coefficient of $i^r/r!$ in the polynomial-like function
$i\mapsto \lg_A(A/q^i)$. The Hilbert-Samuel multiplicity of $A$
itself is defined to be $\fe(m)$. If $A$ is Cohen-Macaulay and $q$
is a parameter ideal (that is, an ideal generated by a system of
parameters), we have $\fe(q)=\lg(A/q)$. If moreover the residue
field is infinite, there exists a parameter ideal $q$
such that $\fe(q)=\fe(m)$ (see exercise 4.5.14 in \cite{BH}).

\begin{lemm} \label{choice_of_parameter_ideal}
Let $k$ be a separably closed field and $X=\Spec(A)$ an affine
scheme of finite type over $k$. Then there exists an integer
$c\ge 1$, a set of Cohen-Macaulay closed points $\cM\subset X$,
and for all points $x\in \cM$ a parameter ideal $q_x\subset\cO_{X,x}$
satisfying $\dim_k(\cO_{X,x}/q_x)\le c$, such that
$$
\bigcap_{x\in\cM}\, q_x'=0
$$
where $q_x'$ is the preimage of $q_x$ in $A$.
\end{lemm}

\begin{proo}
Let $0=I_1\cap\dots\cap I_r$ be a primary decomposition of the
sub-$A$-module $0\subset A$, where $I_j$ is a $P_j$-primary
ideal, $P_j=\sqrt{I_j}$.
For each $1\le j\le r$ let $e_j$ be such that
$(P_j)^{e_j}\subset I_j$. The closed subscheme $Z_j$ defined by the
ideal $P_j$ is a variety, in particular it is reduced.
On one hand, by classical properties of schemes of finite type over
a field, there is a dense open set $U_j\subset Z_j$ of points that
are regular in $Z_j$ and Cohen-Macaulay in $X$. On the other hand,
let $k^\alg$ be an algebraic closure of $k$, and let $S_j$ be the
smooth locus of the reduced subscheme of $Z_j\otimes_k k^\alg$.
It is defined over a finite purely inseparable extension $\ell_j/k$,
whose degree we call $\gamma_j$. Hence there is a smooth
$\ell_j$-scheme $V_j$ whose pullback to $k^\alg$ is $S_j$.
Since $\ell_j$ is separably closed, the set of $\ell_j$-rational
points of $V_j$ is dense. Therefore, the set
$\cM_j=U_j\cap V_j(\ell_j)$ is dense in $Z_j$. By theorem~\ref{thm_EH}
applied with $\cN=\cM_j$ and $M=A/I_j$, we have
$$
\bigcap_{x\in\cM_j}\,m^{e_j}A \subset I_j
$$
where $m$ denotes the maximal ideal of $A$ corresponding to
the point $x$. We call $e=\max(e_j)$, $\cM=\cup\,\cM_j$,
$\gamma=\max(\gamma_j)$. Then, for all $x\in \cM$, we have
$[k(x):k]\le \gamma$, and
$$
\bigcap_{x\in\cM}\, m^eA \ \subset \ I_1\cap\dots\cap I_r=0 \ .
$$
We now choose suitable parameter ideals $q_x$. For
$x\in X$ we let $\fe(x)$ denote the Hilbert-Samuel multiplicity
of the local ring at $x$. This is an upper-semicontinuous
function, hence it is bounded on $X$ by some constant $\alpha$. By the
remarks preceding the lemma, for each Cohen-Macaulay closed point
$x\in X$, we can find a parameter ideal $q=(r_1,\dots,r_s)$ with
$\fe(q)=\fe(x)$, where $s=\dim(\cO_{X,x})\le n=\dim(X)$. Now
$q_x:=((r_1)^e,\dots,(r_s)^e)$ is again a parameter ideal, with
$q_x\subset m^e$. It follows from the above that if $q_x'$ denotes
the preimage of $q_x$ in $A$ then
$$
\bigcap_{x\in\cM}\, q_x'=0 \ .
$$
Furthermore one sees readily that if $\beta=s(e-1)+1$ then
$q^\beta\subset q_x$. Thus,
$$
\lg(\cO_{X,x}/q_x)=\fe(q_x)\le \fe(q^\beta)=\beta^s\fe(q)\le \beta^s\alpha\;.
$$
Finally, since the degree of the residue fields of points
$x\in\cM$ is bounded by $\gamma$, we have
$$
\dim_k(\cO_{X,x}/q_x)=[k(x):k]\,\lg(\cO_{X,x}/q_x)\le
\gamma\beta^s\alpha\le\gamma(n(e-1)+1)^n\alpha\ .
$$
If we set $c:=\gamma(n(e-1)+1)^n\alpha$, we have proven all the
assertions of the lemma.
\end{proo}

\begin{theo} \label{density}
Consider one of the two situations~:
\begin{trivlist}
\itemn{\I} $R$ is henselian and $X$ is an
$R$-scheme locally of finite type, flat and pure.
\itemn{\II} $R$ is complete and $X$ is a flat formal $R$-scheme
locally of finite type.
\end{trivlist}
Then the family $\cF(X)$ of all closed (formal) subschemes
$Z_\lambda\subset X$ finite flat over $R$ is $R$-universally
schematically dense.
\end{theo}

\begin{proo}
We start with case (\I). We first assume that $R$ is strictly henselian.
By lemma~\ref{msep}, there is an open neighbourhood $U$ of the special
fibre of $X$ that is covered by open affine subschemes with
function ring separated for the $\pi$-adic topology.
Lemma~\ref{neighbourhood_sch_dense} implies that $U$ is
$R$-universally schematically dense in $X$. Therefore we may replace
$X$ by $U$ and hence assume that $X$ is covered by open affine
subschemes with function ring separated for the $\pi$-adic topology.
Since the result is local on $X$ we may finally assume that
$X$ is affine, with function ring $A$ of finite type over $R$,
separated for the $\pi$-adic topology (and in fact free, by
remark~\ref{rema_pure_affine_case}).

By lemma~\ref{choice_of_parameter_ideal}, there exists a constant
$c\ge 1$, a set of Cohen-Macaulay closed points $\cM\subset X_k$,
and parameter ideals $q_x\subset\cO_{X_k,x}$ satisfying
$\dim_k(\cO_{X_k,x}/q_x)\le c$ and such that the ideals
$q_x'=q_x\cap A_k$ have zero intersection.
We let $\{Z_\lambda^c\}$, $\lambda\in L^c$, denote the family of all closed
subschemes of $X$, finite flat over $R$, of degree less than $c$,
and we write $f^c:\amalg\,Z_\lambda^c\to X$ for the canonical morphism.

The ideal $q_x$ is generated by a regular sequence $r=(r_1,\dots,r_s)$,
where $s=\dim(\cO_{X_k,x})$. Let $\tilde r$ be a sequence obtained by
lifting the $r_i$ in $\cO_{X,x}$ and let $Y=\Spec(\cO_{X,x}/(\tilde r))$.
As $r$ is a regular sequence, it follows that $Y$ is flat over $R$.
Furthermore $Y_k$ is artinian, hence $Y$ is quasi-finite over~$R$. Since
$R$ is henselian, $Y$ is in fact finite over $R$. Thus $Y\to X$ is a
proper monomorphism, hence a closed immersion. So $Y$ is one of the
schemes $Z_\lambda^c$.

Since the $k$-algebras of functions of $Z_{\lambda,k}^c$ are free of
rank less than $c$, the Cayley-Hamilton theorem implies that
in the terminology of definition~\ref {torsion_bounded},
the torsion in $(f^c_k)_*\cO_{Z_{\lambda,k}^c}$ is bounded uniformly
in $\lambda$ (in a strong form, since the bound $c$ is independent of the local
sections $t$). As the intersection of the ideals $q_x'=q_x\cap A_k$
is zero, lemma~\ref{criterion_bounded_torsion} applies and proves
that $f^c_k$ is schematically dominant. Moreover, the $R$-algebras of
functions of $Z_\lambda^c$ are free of rank less than~$c$, so the
argument used above works again and by
proposition~\ref{cocartesian_diagram_affine_case}
we get that $f^c$ and $f^c_K$ are schematically dominant. Applying
\cite{EGA}~IV, 11.10.9, it follows that $f^c$ is $R$-universally
schematically dominant. A fortiori, the family $\cF(X)$
is $R$-universally schematically dense.

It remains to treat the case of a general henselian discrete
valuation ring $R$. Let $R^\sh$ be a strict henselization, and
$X^\sh=X\otimes_R R^\sh$. By the preceding discussion we know that
$\cF(X^\sh)$ is universally schematically dense in $X^\sh$. Since
$R^\sh$ is an integral extension of $R$, the canonical morphism
$j\colon X^\sh\to X$ is integral. Thus the schematic image of any
finite $R^\sh$-flat closed subscheme $Z^\sh\subset X^\sh$ is an
$R$-flat closed subscheme $Z$ of $X$, integral over $R$, hence a
finite flat $R$-scheme. This proves that the family $\{j^{-1}(Z)\}$,
with $Z\in\cF(X)$, is a cofinal subfamily of $\cF(X^\sh)$, thus
it is universally schematically dense in $X^\sh$. By faithfully
flat descent (\cite{EGA}~IV.11.10.5), so is $\cF(X)$ in $X$.

In case (\II), we follow the same strategy of proof. We start with the
case where $R$ is strictly henselian. We reduce to the formal affine
case $X=\Spf(A)$, with $A$ topologically of finite type over $R$. Such
an $A$ is automatically separated for the $\pi$-adic topology. Then
we consider the family $\{Z_\lambda^c\}$ of all closed formal
subschemes of $X$, finite flat over $R$, of degree less than $c$.
We apply lemma~\ref{choice_of_parameter_ideal} again, and as before,
for each Cohen-Macaulay closed point $x$ in $\cM\subset X_k$, we can
realize the subscheme defined by the parameter ideal
$q_x\subset\cO_{X_k,x}$ as the special fibre of some $Z_\lambda^c$.
Then we use proposition~\ref{cocartesian_diagram_affine_formal_case}
to get that $f^c$ and $f^c_K$ are schematically dense. It makes no
difficulty to adapt \cite{EGA}~IV, 11.10.9 to formal schemes and
conclude that $f^c$ and a fortiori $\cF(X)$ is $R$-universally
schematically dense. Also the argument from \cite{EGA} to descend
from the strict henselization to $R$ is easily adapted.
\end{proo}

\begin{prop} \label{cocartesian_diagram_general_case}
Let $X$ be an $R$-scheme locally of finite type and flat.
Let $\{Z_\lambda\}$ be the family of all closed subschemes
of $X$ finite flat over $R$, and assume that the family
$\{Z_{\lambda,k}\}$ is schematically dense in $X_k$ (e.g.
$R$ is henselian and $X$ is pure, by theorem~\ref{density}).
Then for all separated $R$-schemes $Y$ and all diagrams in
solid arrows
$$
\xymatrix{\amalg Z_{\lambda,K} \ar[r] \ar[d] & X_K \ar[d] \ar[rdd]^\beta & \\
\amalg Z_\lambda \ar[r] \ar[rrd]_\alpha & X \ar@{.>}[rd] & \\
& & Y}
$$
there exists a unique morphism $g:X\to Y$ making the full
diagram commutative.
\end{prop}

\begin{proo}
In fact flatness of $X$ follows from the other assumptions,
by proposition~\ref{cocartesian_diagram_affine_case}. Let
$f:\amalg\,Z_\lambda\to X$, $\alpha:\amalg\,Z_\lambda\to Y$ and
$\beta:X_K\to Y$ be the maps in the diagramme. By the same
argument as in the proof of
proposition~\ref{cocartesian_diagram_affine_case}, the map
$g:X\to Y$ is unique if it exists. Thus we may define it
locally on $X$ and glue. It is enough to define $g$ in a
neighbourhood of all closed points $x\in X_k$. We know from
lemma~\ref{closed_points_attained} that $f$ is surjective on
closed points of the special fibre. So the given point $x$ is
equal to $f(z)$ for some $\lambda$ and $z\in Z_\lambda$. Let
$y=\alpha(x)$, let $V=\Spec(C)$ be an open affine neighbourhood
of $y$ in $Y$, and let $U$ be an open subscheme of $X$
containing~$x$. We will prove that $x$ does not belong to the
closure in $X$ of $X_K\setminus \beta^{-1}(V)$. Indeed,
otherwise there is a point $\eta\in X_K\setminus \beta^{-1}(V)$
such that $x\in W:=\overline{\{\eta\}}$. Thanks to
lemma~\ref{closed_points_attained} applied to~$W$, we may
replace $\eta$ by a closed point of $W_K$ and hence assume
that $\eta$ is closed in $X_K$. In this case~$W$ is one of the
$Z_\lambda$, so it makes sense to speak about the images of $x$
and $\eta$ under $\alpha$. Then,
$$
x\in \overline{\{\eta\}} \quad\mbox{implies that} \quad
y=\alpha(x)\in \overline{\{\alpha(\eta)\}}=\overline{\{\beta(\eta)\}}
$$
and this is a contradiction with the fact that
$\beta(\eta)\not\in V$. Therefore, we may shrink $U$ and
assume that $U_K\subset \beta^{-1}(V)$. Then by
lemma~\ref{msep} we may shrink $U$ further to the
spectrum of a ring~$A$ separated for the $\pi$-adic topology.
Therefore, we reduce to $X=\Spec(A)$ and $Y=\Spec(C)$, and
proposition~\ref{cocartesian_diagram_affine_case} applies.
\end{proo}

\begin{prop} \label{cocartesian_diagram_general_formal_case}
Assume that $R$ is complete. Let $X$ be a flat formal $R$-scheme of
finite type. Let $\{Z_\lambda\}$ be the
family of all closed formal subschemes of $X$ finite flat
over $R$ and $f:\amalg\,Z_\lambda\to X$ the canonical map.
Then the analogue of the amalgamated sum property of
proposition~\ref{cocartesian_diagram_general_case} holds,
if we understand a morphism from a rigid analytic
$K$-space $Z$ to a formal $R$-scheme $Y$ to be a morphism
$Z\to Y_\rig$. More precisely, given
\begin{itemize}
\item a separated formal $R$-scheme $Y$,
\item a morphism of formal $R$-schemes
$\alpha:\amalg\,Z_\lambda\to Y$,
\item a morphism of rigid spaces $\beta:X_\rig\to Y_\rig$
\end{itemize}
such that $\alpha_\rig=\beta\circ f_\rig$, there exists a unique
morphism $g:X\to Y$ such that $g_\rig=\beta$ and $g\circ f=\alpha$.
\end{prop}

\begin{proo}
The proof of \ref{cocartesian_diagram_general_case} works again in
this setting, with some adaptations which we now sketch.
If $g,g'\colon X\to Y$ are two solutions to the problem,
then in particular $g_\rig=g'_\rig$. By Raynaud's theorem
(\cite{BL1}, th. 4.1) there exists an
admissible formal blowing-up $s:X'\to X$ such that $g\circ s=g'\circ
s$. Since $s$ is schematically dominant and $Y$ is separated, we get
$g=g'$. Because of this unicity statement, as far as existence is
concerned, we may define $g$ locally on $X$ and glue.
Also we need to know that $\cF(X)_k$ is schematically dense, which
is granted by theorem~\ref{density}. Then, by the same method
as above, we reduce to the affine case $X=\Spf(A)$ and $Y=\Spf(C)$.
Now the arguments of the proof of point (3) in
proposition~\ref{cocartesian_diagram_affine_case} carry on.
\end{proo}

\section{Schematic images inside $\Hom$ and $\Aut$ functors}
\label{Sch_images_inside_Hom_and_Aut}




Throughout this section, we fix a discrete valuation ring $(R,K,k,\pi)$.
We first recall the definition of schematic closures and images for
fppf sheaves over a discrete valuation ring $R$. After a brief
discussion of kernels, we prove the main theorems of the paper
on representability of schematic images.

\subsection{Definitions} \label{def_and_example}

Recall that in the context of schemes, if $f\colon W\to X$ is a
morphism of $R$-schemes such that $f_*\cO_W$ is quasi-coherent,
there exists a smallest closed subscheme $X'\subset X$ such that
$f$ factors through~$X'$. We call it the {\em schematic image}
of $f$. If is equivalent to say that the schematic image of $f$
is~$X$, or that $f$ is schematically dominant.

If $W$ is a closed subscheme of the generic fibre of $X$ and
$f$ is the canonical immersion, then the schematic image is
called the {\em schematic closure} of $W$ in $X$. It is the
unique closed subscheme of $X$ which is flat over $R$ and whose
generic fibre is $W$ (\cite{EGA}~IV.2.8.5).

These definitions may be adapted to morphisms of sheaves as follows
(see \cite{Ra})~:

\begin{defis} \label{defi_AS}
Let $F$ be an fppf sheaf over the category of $R$-schemes.
\begin{trivlist}
\itemn{1} Let $G$ be a subsheaf of the generic fibre $F_K$. Then the
{\em schematic closure of $G$ in $F$} is the fppf sheaf $G'$
associated to the presheaf $G^\flat$ defined as follows. Given an
$R$-scheme $T$, $G^\flat(T)$ is the set of all morphisms
$f\colon T\to F$ such that there exists a factorization
$$
\xymatrix{T\ar[r]
  \ar[rd]_f & T' \ar[d]^g \\ & F}
$$
with $T'$ flat over $R$ and $g(T'_K)\subset G$.
\itemn{2} We say that $F$ is {\em flat over $R$} if it is
equal to the schematic closure of its generic fibre.
\itemn{3} Let $h\colon H\to F$ be a morphism of fppf sheaves over
$R$, with $H$ flat. Let $G$ be the image sheaf of $h_K\colon H_K\to
F_K$. Then the {\em schematic image} of $H$ in $F$ is defined to be the
schematic closure of $G$ inside $F$.
\end{trivlist}
\end{defis}

The following properties are formal consequences of the definitions.
The formation of the schematic closure commutes with flat extensions
of discrete valuation rings. Let $F_1,F_2$ be sheaves over the
category of $R$-schemes. Let
$G_1\subset F_{1,K}$, $G_2\subset F_{2,K}$ be subsheaves, and let
$G_1',G_2'$ be the schematic closures. For a morphism of sheaves
$\alpha\colon F_1\to F_2$ such that $\alpha(G_1)\subset G_2$, we have
$\alpha(G_1')\subset G_2'$. As a consequence, the schematic closure of
$G$ in $F$ is the only subsheaf of $F$
which is flat over~$R$ and has generic fibre equal to $G$. Finally the
formation of the schematic closure commutes with products~; it follows
that if $F$ is a group (resp. monoid) sheaf i.e. a group (resp. monoid)
object in the category of fppf sheaves, and $G$ is a subgroup
(resp. submonoid) sheaf of $F_K$, then the schematic closure $G'$ is a
subgroup (resp. submonoid) sheaf of~$F$.

In general, even if $F$ is representable by a scheme, one needs
rather strong conditions on the monomorphism $G\to F_K$ if one
wants representability of the schematic closure $G'$. As we
recalled above, one pleasant case is when $G\to F_K$ is a closed
immersion~; then $G'\to F$ is also a closed immersion. However,
we will see now that already in the case of an open immersion,
the schematic closure is not representable by a scheme in general.

\begin{lemm}
Let $X$ be an $R$-scheme, $U_K\subset X_K$ the complement of a
Cartier divisor. Then the schematic closure of $U_K$ in
$X$ is representable by an inductive limit of affine $X$-schemes.
\end{lemm}

\begin{proo}
We first construct $U'$. Fix an integer $n\ge 0$. For each open affine
$V=\Spec(A)$ in $X$, me may choose an equation $f\in A$ for $X_K\setminus
U_K$. Define $U_{V,n}$ to be the spectrum of the ring
$$
\left(\frac{A[x_n]}{x_nf-\pi^n}\right)_0
$$
where the subscript $0$ means the quotient by the $\pi$-torsion
ideal $(0:\pi^\infty)$. There are maps $U_{V,n}\to U_{V,n+1}$ given by
$x_{n+1}\mapsto \pi x_n$, and we define $U'_V$ to be the limit of the
schemes $U_{V,n}$. This construction glues over all $V$ to give an
$X$-scheme $U'$. It is not hard to see that this is independent of the
choice of local equations $f$, up to isomorphism. Finally we check
that $U'$ is the desired schematic closure. Let $g\colon T\to X$ be a
morphism with $T$ flat over $R$ and $g(T_K)\subset U_K$. Let
$V=\Spec(A)$ be an open affine in $X$ and $W=\Spec(B)$ open affine in
$T$, with $g(W)\subset V$~; let $f\in A$ be an equation for
$X_K\setminus U_K$. Then we have a morphism of rings $\varphi\colon
A\to B$ such that $\varphi(f)$ is invertible in $B_K$, i.e. there
exists $n\ge 0$ and $t\in B$ such that
$\varphi(f)t=\pi^n$. Furthermore, since $B$ is $R$-flat, $t$ is
uniquely determined, as well as the morphism of $A$-algebras
$$
\left(\frac{A[x_n]}{x_nf-\pi^n}\right)_0\to B
$$
given by $x_n\mapsto t$. These morphisms glue to a unique map $T\to U'$.
\end{proo}

\subsection{Kernels} \label{kernels}

Let $S$ be a base scheme and let $\Gamma,X,Y$ be schemes over $S$. We
consider a morphism of $S$-schemes $\varphi\colon \Gamma\times_S X\to Y$,
which we view as an action of $\Gamma$ on $X$ with values in
$Y$. Equivalently, we have a morphism of functors $\varphi'\colon
\Gamma\to\Hom_S(X,Y)$. We say that {\em $\Gamma$ acts faithfully on $X$},
or that {\em $\varphi$ is faithful}, if $\varphi'$ is a monomorphism.
We can relate this to the morphism~:
$$
\varphi'\times\varphi'\colon \Gamma\times_S \Gamma\to \Hom_S(X,Y)\times_S
\Hom_S(X,Y) \ .
$$

\begin{defi}
The {\em kernel} of $\varphi$ is the
preimage of the diagonal of $\Hom_S(X,Y)$ under the morphism
$\varphi'\times\varphi'$. It is denoted $\ker(\varphi)$.
\end{defi}

Obviously, it is equivalent to say that $\Gamma$ acts faithfully on $X$,
or that the natural monomorphism $\Delta\to\ker(\varphi)$ is
an isomorphism, $\Delta\subset \Gamma\times_S\Gamma$ being the diagonal
of $\Gamma$. When this holds, we shall also say by abuse
that {\em $\ker(\varphi)$ is trivial}.

If $X=Y$ and $G$ is a group scheme acting on $X$,
the relation between the kernel we have just
defined and the usual kernel $H:=(\varphi')^{-1}(\Id_X)$ is
given by the isomorphism $G\times H\to \ker(\varphi)$
taking $(g,h)$ to $(g,gh)$. We use the notation
$\ker(\varphi)$ in both situations, because the context will
never allow confusions.

The lemma below collects some cases where one knows that the kernel
is representable by a closed subscheme of $\Gamma\times_S\Gamma$.
One of this case involves essentially free morphisms of schemes,
a notion which can be slightly (and fruitfully) generalized to
essentially semireflexive (see \cite{SGA3}, Expos{\'e}~VIII, \S~6
and \cite{To2}, \S~1). Recall that a module $M$ over a ring $A$
is called
{\em semireflexive} if the natural morphism $M\to M^{\vee\vee}$
to the linear bidual is injective. It is equivalent to say that
$M$ can be embedded into a product module $A^I$, for some set $I$.
A morphism of schemes $X\to S$ is called {\em essentially free}
(resp. {\em essentially semireflexive}) over $S$, if there
exists a covering of $S$ by open affine schemes $S_i$, for all $i$
an affine scheme $S'_i$ faithfully flat over $S_i$, and a covering
of $X'_i=X\times_S S'_i$ by open affine schemes $X'_{i,j}$, such that
for all $i,j$ the function ring of $X'_{i,j}$ is a free (resp.
semireflexive) module over the function ring of $S'_i$.
It is clear that an essentially free morphism is essentially
semireflexive.

\begin{lemm} \label{lemma_representability_of_kernel}
Let $X\to S$ be flat and $Y\to S$ separated. Then
$\ker(\varphi)\to \Gamma\times_S \Gamma$ is a closed immersion
in any of the following cases~:
\begin{trivlist}
\itemn{i} $X\to S$ is essentially semireflexive,
\itemn{ii} $S$ is regular noetherian of dimension $1$
and $X\to S$ is locally of finite type, flat and pure,
\itemn{iii} $X\to S$ is proper and $X,Y$ are locally of finite
presentation over $S$.
\end{trivlist}
\end{lemm}

We see that under one of these three conditions, faithfulness of
$\varphi$ implies separation of $\Gamma$. We remark also that it
is not hard to see that if $X\to S$ is flat and $Y\to S$ is separated,
then $\ker(\varphi)\to \Gamma\times_S \Gamma$ satisfies the
valuative criterion of properness. What is more difficult is to
check that it is of finite type.

\begin{proo}
For case (i) we refer to \cite{SGA3}, Expos{\'e}~VIII, \S~6
and \cite{To2}, lemma 1.16.
In case (iii), the functor $\Hom_S(X,Y)$ is a separated
algebraic space, by Artin's theorems, so the result is clear.
It remains to consider case~(ii). We may assume that $S$ is the spectrum
of a henselian discrete valuation ring $R$. By lemma~\ref{msep},
there is an open neighbourhood $U$ of the special fibre of $X$
that is covered by open affine subschemes $U_i$ with function
ring $A_i$ separated for the $\pi$-adic topology.
Besides, $U_i$ is pure over $R$ and $A_i$ is a free $R$-module, by
\ref{rema_pure_affine_case}. It follows that $U$ is essentially free
over $S$, hence the kernel $N_U:=\ker(\Gamma\times_S U\to Y)$
is a closed subscheme of $\Gamma\times_S\Gamma$ by case (i).
Consider the map induced by the action~:
$$
\psi:N_U\times_S X\to Y\times_S Y
$$
given on the points by
$$
(\gamma_1,\gamma_2,x)\mapsto (\varphi(\gamma_1)(x),\varphi(\gamma_2)(x)) \ .
$$
By definition, the restriction of $\psi$ to $N_U\times_S U$
factors through the diagonal of $Y$. Since $U$ is $R$-universally
schematically dense in $X$ (lemma~\ref{neighbourhood_sch_dense}),
then $N_U\times_S U$ is schematically dense in $N_U\times_S X$.
Thus $\psi$ factors through the diagonal, that is,
$N_U\to \Gamma\times_S\Gamma$ factors through the kernel
$N:=\ker(\Gamma\times_S X\to Y)$. This gives an
inverse for the obvious morphism
$N\to N_U$ and proves that $N\simeq N_U$. In particular, $N$ is a
closed subscheme of $\Gamma\times_S\Gamma$, as claimed.
\end{proo}

\begin{lemm} \label{proper_gp_good_action}
Let $X$, $Y$, $\Gamma$ be $R$-schemes.
Consider one of the two situations~:
\begin{trivlist}
\itemn{1} $R$ is henselian, $X$ is locally of finite type, flat
and pure, $Y$ is separated, $\Gamma$ is noetherian.
\itemn{2} $X$ is affine and the family of its closed subschemes
finite flat over $R$ is universally affinely dominant
(definition~\ref{defis_dominance}), $Y$ is affine, $\Gamma$ is noetherian.
\end{trivlist}
Consider an action $\varphi:\Gamma\times X\to Y$ faithful on the
generic fibre. Then there exists a finite $R$-flat closed subscheme
$Z\subset X$ such that the induced action $\Gamma\times Z\to Y$ has
the same kernel as $\varphi$.
\end{lemm}

\begin{proo}
Note that in case (2), the scheme $X$ is semireflexive, so that
in both cases the kernels are closed subschemes of $\Gamma\times\Gamma$
by lemma~\ref{lemma_representability_of_kernel}.
Let $N=\ker(\varphi)$. Let
$Z_1\subset X$ be a finite flat closed subscheme and let $N_1$ be the
kernel of the restricted action $\Gamma\times Z_1\to Y$. If $N_1\ne
N$, there exists $Z_2\subset X$ with $Z_1\subset Z_2$ such that
$N_1\supsetneq N_2$. For, otherwise $N_1$ would act trivially on
all the finite flat closed subschemes $Z\supset Z_1$, which are
universally schematically dense in $X$ (by theorem~\ref{density}
in case (1)), hence $N_1$ would act trivially on $X$, a contradiction.
For $s\ge 1$, as long as $N_s\ne N$, we iterate this process and
obtain a sequence
$N_1\supsetneq N_2\supsetneq N_3\supsetneq\dots$
Since $\Gamma\times \Gamma$ is noetherian, for some $s$
we obtain that $N_s=N$. We can choose $Z=Z_s$.
\end{proo}

\subsection{Representability of schematic images}
\label{main_result}

We now come to the main results of this paper.

\begin{lemm} \label{representability_1}
Let $R$ be a discrete valuation ring. Let $X,Y$ be $R$-schemes
locally of finite type, with $X$ flat and pure and $Y$ separated.
Consider a finite flat $R$-scheme $\Gamma$ and an action
$\varphi:\Gamma\times X\to Y$ faithful on
the generic fibre. Then the schematic image of $\Gamma$ in
$\Hom_R(X,Y)$ is representable by a finite flat $R$-scheme
$\Gamma'$.
\end{lemm}

We stress again that $\Hom_R(X,Y)$ is far from being
representable, in general.

\begin{proo}
We start with the case where $R$ is henselian.
By lemma~\ref{proper_gp_good_action}, there is a finite $R$-flat
closed subscheme $Z_0\subset X$ such that $\Gamma_K$ acts faithfully
on $Z_{0,K}$. Let $\{Z_\lambda\}_{\lambda\in L}$ be the family
of all finite $R$-flat closed subschemes of $X$ containing $Z_0$. This
family carries the filtering order by inclusion of subschemes~:
$\lambda\le \mu$ if and only if $Z_\lambda\subset Z_\mu$. Since
$Z_\lambda$ is finite flat over $R$, the functor $\Hom_R(Z_\lambda,Y)$ is
representable by a scheme. Moreover, since $Z_\lambda\supset Z_0$
and $\Gamma$ is finite, the map $\Gamma_K\to\Hom_K(Z_{\lambda,K},Y_K)$
is a closed immersion. For each $\lambda$ we define $\Gamma'_\lambda$
to be the schematic image of the map $\Gamma\to\Hom_R(Z_\lambda,Y)$.
If $\lambda\le\mu$ in $L$, there is a restriction morphism
$\Hom_R(Z_\mu,Y)\to\Hom_R(Z_\lambda,Y)$ and taking schematic closures
gives maps $\Gamma'_\mu\to \Gamma'_\lambda$. Let $\Gamma'$ be the
filtering projective limit of the system $\{\Gamma'_\lambda\}$. This
is an affine, flat, integral $R$-scheme~; it is dominated by
$\Gamma$ hence finite over $R$. Applying
proposition~\ref{cocartesian_diagram_general_case} to the diagram
$$
\xymatrix{\amalg\; \Gamma'_K\times Z_{\lambda,K} \ar[d] \ar[r]^\simeq &
\Gamma'_K\times \amalg\; Z_{\lambda,K}
  \ar[r] & \Gamma'_K\times X_K \ar[d] \ar@/^1pc/[rdd] & \\
\amalg\; \Gamma'\times Z_\lambda \ar[r]^\simeq \ar@/_1pc/[rrrd] &
\Gamma'\times \amalg\; Z_\lambda \ar[r] &
\Gamma'\times X \ar@{.>}[rd] & \\
& & & Y}
$$
we obtain an action of $\Gamma'$ on $X$ with values in $Y$. This action is
clearly universally faithful, because the morphism $\amalg\,\Gamma'\times_R
Z_\lambda\to \Gamma'\times_R X$ is universally schematically dominant (apply
theorem~\ref{density} to $X$ and pull back to $\Gamma'\times_R X$). So
$\Gamma'$ has the characterizing properties of the schematic closure of
$\Gamma$ in $\Hom_R(X,Y)$, and this proves the theorem.

If $R$ is an arbitrary discrete valuation ring,
let $R^\hh$ be a henselization of $R$. By the preceding discussion,
$\Gamma\otimes_R R^\hh$ is representable by a finite flat
$R^\hh$-scheme. So by descent using \cite{BLR}~6.2/D.4, $\Gamma$ is
representable by a finite flat $R$-scheme.
\end{proo}

There is also a version in the affine case, where one can relax
the assumptions of finite type. For example, it applies to rings
arising from the completion of smooth $R$-schemes along a section.

\begin{lemm} \label{representability_2}
Let $X$ be an affine flat $R$-scheme such that the family of
its closed subschemes finite flat over $R$ is universally
affinely dominant (definition~\ref{defis_dominance}). Let $Y$ be
an affine $R$-scheme and $\Gamma$ an $R$-scheme locally of finite
type, flat and pure. Consider an
action $\varphi:\Gamma\times X\to Y$ faithful on the generic
fibre. Then the schematic image of $\Gamma$ in $\Hom_R(X,Y)$ is
representable by a flat $R$-scheme which is affine if $\Gamma$
is, and finite if $\Gamma$ is.
\end{lemm}

\begin{proo}
Observe that the assumptions imply that $X$ is semireflexive
over $R$, therefore kernels of actions are representable by closed
subschemes, by lemma~\ref{lemma_representability_of_kernel}.
Let $X=\Spec(A)$ and $Z_\lambda=\Spec(B_\lambda)$,
$\lambda\in L$, be the family of the finite flat closed
subschemes of $X$, and let $B=\Pi\,B_\lambda$.
The proof goes in three steps.

\medskip

\no {\em First step~: $\Gamma$ is finite.} In this case we follow
the proof of lemma~\ref{representability_1}. All references to
theorem~\ref{density} are replaced by the assumption made on $X$.
The reference to proposition~\ref{cocartesian_diagram_general_case}
is replaced by a reference to
proposition~\ref{cocartesian_diagram_affine_case}.
The conclusion is that the schematic image is representable by
a finite flat $R$-scheme $\Gamma'$.

\medskip

\no {\em Second step~: $\Gamma$ is affine.} Let $\Gamma=\Spec(C)$ and
call $\Delta_\mu=\Spec(D_\mu)$, $\mu\in M$ the family of all finite $R$-flat
closed subschemes of $\Gamma$. By the first step, for all $\mu$
the schematic image of $\Delta_\mu$ in $\Hom_R(X,Y)$ is representable
by a finite flat $R$-scheme $\Delta'_\mu=\Spec(D'_\mu)$. Let
$D=\Pi\,D_\mu$, $D'=\Pi\,D'_\mu$. We have injective ring
homomorphisms $C\into D$ and $D'\into D$. Let $C'$ be the
intersection of $C$ and $D'$ inside $D$, and $\Gamma'=\Spec(C')$.
We claim that $\{\Delta'_\mu\}_{\mu\in M}$ is the family of all
finite flat closed subschemes
of $\Gamma'$. Indeed, it is easy to see that $C'\to D'_\mu$ is
surjective, i.e. $\Delta'_\mu$ is a finite flat closed subscheme of
$\Gamma'$. Moreover, for each finite flat closed subscheme
$T'\subset \Gamma'$, we can consider $T'_K$ as a closed subscheme
of $\Gamma_K$, we set $\Delta_\mu$ equal to the schematic closure of
$T'_K$ in $\Gamma$, then obviously $T'=\Delta'_\mu$. Now we prove
that $\Gamma'$ acts on $X$. For this, note that $\coker(C'\to D')$
injects into $\coker(C\to D)$ and hence is $R$-flat. It follows
from lemma~\ref{equivalent_conditions_univ_injective} that the
family of finite flat closed subschemes of $\Gamma'$ is universally
affinely dominant.
Then the affine scheme $\Gamma'\times X$ has a family of
finite flat subschemes $\Delta'_\mu\times Z_\lambda$ which
is universally affinely dominant. Using
proposition~\ref{cocartesian_diagram_affine_case}, one obtains
an action $\Gamma'\times X\to Y$. It is clear that this action
has trivial kernel, hence $\Gamma'$ is the schematic image of
$\Gamma$.

\medskip

\no {\em Third step~: $\Gamma$ is arbitrary.} By
lemma~\ref{purity_versus_adic_separation} and lemma~\ref{msep},
there is an open neighbourhood of the special fibre of $\Gamma$
that is covered by pure open
affine subschemes $U_i$. For each $i$, by the second step the
schematic image of $U_i$ is representable by an affine flat
$R$-scheme $U'_i$. By unicity of the schematic image, the
formation of $U'_i$ is compatible with localisation, so that
the various $U'_i$ glue to give a flat $R$-scheme $U'$. Since
$U'_K\simeq U_K$ we can glue $U'$ and $\Gamma_K$ along their
intersection to get a flat $R$-scheme $\Gamma'$. It is clear
that this is the schematic image of $\Gamma$.
\end{proo}

In the sequel, we examine the most interesting case of images
of {\em groups} acting on schemes by group homomorphisms. We
introduce some terminology.

\begin{defi}
If an $R$-group scheme $G$ acts on an $R$-scheme $X$ in such a
way that the action on the generic fibre is faithful, then the
schematic image of $G$ in $\Aut_R(X)$ is called the
{\em effective model} of $G$ for its action on $X$.
\end{defi}

\begin{theo} \label{representability_group}
Let $X$ be an affine flat $R$-scheme whose closed subschemes
finite flat over $R$ form a universally affinely dominant family.
Let $G$ be an $R$-group scheme locally of finite type, flat and pure,
acting on $X$, faithfully on the generic fibre.
Then the effective model $G'$ of the action is representable
by a flat $R$-group scheme. If $G$ is quasi-compact, or affine,
or finite, then $G'$ has the same property.
\end{theo}

\begin{proo}
Let $G''$ be the schematic image of $G$ inside $\Hom_R(X,X)$.
By the previous lemma $G''$ is representable by a flat
$R$-scheme. Since $\Aut_R(X)$ is an open subfunctor of $\Hom_R(X,X)$,
the preimage of $G''$ in $\Aut_R(X)$ is flat over $R$ and hence
is the schematic image $G'$. It follows from the general remarks of
subsection~\ref{def_and_example} that $G'$ is a sub-$R$-group
scheme of $\Aut_R(X)$.

If $G$ is quasi-compact, let $(U'_i)_{i\in I}$ be an open cover
of $G'$. Let $U_i$ be the preimage of $U'_i$ in~$G$. By assumption,
a finite number of open sets $U_1,\dots,U_n$ cover $G$. The scheme
$G'$ is covered by the schematic images of $U_1,\dots,U_n$ which
are none other than $U'_1,\dots,U'_n$. It follows that $G'$ is
quasi-compact.

If $G$ is affine, then $G''$ is affine by
lemma~\ref{representability_2}, hence $G'$ is quasi-affine.
Let $H$ be the affine hull of $G'$. This is a
flat group scheme containing $G'$ as an open subgroup. Moreover,
the special fibre $G'_k$ is schematically dense in the special
fibre $H_k$, and since these are $k$-group schemes, we have
in fact $G'_k=H_k$. It follows that $G'=H$ is affine.

If $G$ finite, then $G\to G'$ is surjective and it follows
easily that $G'$ is finite.
\end{proo}

These representability results extend obviously to the case where
$X$ is covered by invariant open affine subschemes satisfying
the relevant assumptions. When $X$ is locally of finite type
but not necessarily affine, it is more difficult to prove
that schematic images are representable. In fact, it is easy to
provide a group scheme $G^c$ which is a candidate to be the image,
but in order to prove that it acts on $X$ using
proposition~\ref{cocartesian_diagram_general_case}, one needs
$G^c$ to be of finite type. This is the major difficulty of our
method. Moreover, it seems that in numerous situations one can
not expect the schematic image $G'$ to be of finite type unless the
kernel of the action of~$G$ is very small. The following two
results give examples of this.

\begin{theo} \label{representability_group_proper_case}
Let $X$ be an $R$-scheme locally of finite type, separated,
flat and pure. Let $G$ be a flat proper $R$-group scheme
acting on $X$, faithfully on the generic fibre. Let $N$
denote the kernel of the action. Then the effective model
$G'$ is representable by a flat group scheme of finite type
if and only if $N_k$ is finite. Moreover, in this case $G'$
is proper.
\end{theo}

\begin{proo}
First, assume that $N_k$ is finite. We adapt the proof of
lemma~\ref{representability_1}. By
lemma~\ref{proper_gp_good_action}, there is a finite $R$-flat
closed subscheme $Z_0\subset X$ such that $G_K$ acts faithfully
on $Z_{0,K}$. Let $G_0$ be the schematic image of $G$ inside
$\Hom_R(Z_0,X)$, which is representable since $Z_0$ is finite.
We claim that the morphism $u\colon G\to G_0$ is finite. Indeed,
on the special fibre $u_k$ factors as the composition of the finite
quotient $G_k\to G_k/N_k$ and the monomorphism $G_k/N_k\to G_{0,k}$
given by the embedding in $\Hom_k(Z_k,X_k)$. It follows that $u$ is
quasi-finite, hence finite since $G$ is proper.

Now let $\{Z_\lambda\}_{\lambda\in L}$ be the family of all finite
$R$-flat closed subschemes of $X$ containing $Z_0$. For each
$\lambda$, let $G''_\lambda$ be the schematic image of the
map $G\to\Hom_R(Z_\lambda,X)$. Since $G\to G''_\lambda\to G_0$
is finite and schematically dominant, then $G\to G''_\lambda$
and $G''_\lambda\to G_0$ are finite schematically dominant also.
Let $G''$ be the filtering projective limit of the system
$\{G''_\lambda\}$. This is a scheme which is finite over $G_0$. Also,
$G\to G''$ is finite, thus $G''$ is of finite type over~$R$ by the
Artin-Tate theorem (see \cite{Ei}, exercise 4.32). Applying
proposition~\ref{cocartesian_diagram_general_case} like in the
proof of theorem~\ref{representability_1}, we obtain an action
of $G''$ on $X$ with values in $X$. Let $G'$ be the preimage of
$G''$ under the inclusion $\Aut_R(X)\subset\Hom_R(X,X)$. This
is the schematic image of $G$ in $\Aut_R(X)$. Since $G\to G'$
is finite, then $G'$ is proper.

Conversely, assume that $G'$ is representable by a flat group
scheme of finite type over $R$. A result of Anantharaman
asserts that a separated morphism $u$ between flat $R$-group schemes
of finite type such that $u_K$ has affine kernel is affine (\cite{An},
chap.~II, prop.~2.3.2). It follows that $G\to G'$ is
affine. Since it is also proper, it is in fact finite. It follows
easily that $N_k$ is finite.
\end{proo}

\begin{rema}
It is a well-known fact that a proper flat group scheme over $R$
is in fact projective. Here is one way to see it. Given a finite
extension $K^*/K$, write $G^*$ for the extension of $G$ to $R^*$,
the integral closure of $R$ in $K^*$. By a result of Raynaud and
Faltings (\cite{PY}, corollary~A.4) there is a finite extension $K^*/K$
such that the normalization morphism $(\tilde G^*)_\red\to (G^*)_\red$
is finite and $(\tilde G^*)_\red$ is smooth. Hence it is the product of
an abelian scheme by an {\'e}tale finite group, hence projective. It
follows that $(G^*)_\red$ is projective, hence also $G^*$ and $G$
itself. Another way to check that $G$ is projective is to reduce to
the connected case. Then $G$ is commutative and one can apply
\cite{An}, chap.~II, prop.~2.2.1.
\end{rema}

\begin{rema}
Under the assumptions of
theorem~\ref{representability_group_proper_case}, it seems
plausible that if $N_k$ is finite, then $G'$ is representable
whether $G$ is proper or not. The only point that needs a
verification is that $u:G\to G_0$ is finite (with the notations
of the proof of the proposition). Even though $u_K$ and $u_k$ are
finite, I was not able to prove this.
\end{rema}

\begin{prop}
Let $X$ be an $R$-scheme locally of finite type, separated, flat
and pure. Let $G$ be a reductive $R$-group scheme acting on $X$,
faithfully on the generic fibre. Assume furthermore that
either $k$ has characteristic $p\ne 2$, or that no normal
subgroup of $G_{\bar K}$ is isomorphic to $\SO_{2n+1}$ for some
$n\ge 1$. Let $N$ denote the kernel of the action. Then the
effective model $G'$ is representable by a flat group scheme
of finite type if and only if $N=1$.
\end{prop}

\begin{proo}
This is in fact a rigidity property of reductive groups. Assume
that $G'$ is representable by a flat group scheme of finite type.
Since $X$ is flat and separated, then $\Aut_R(X)$ is a separated
sheaf. It follows that $G'$ is separated. Then $G'$ is
affine by \cite{An}, chap.~II, prop.~2.3.1. By corollary~1.3
of \cite{PY}, we obtain that $G\to G'$ is a closed
immersion. It follows that $G$ acts faithfully on $X$, in other
words $N=1$. The converse is obvious.
\end{proo}

From this proposition follows that if $G$ is a finite group
scheme of order prime to $p=\ch(k)$ acting on an $R$-scheme
locally of finite type, separated, flat and pure $X$, then $G$
acts faithfully as soon as $G_K$ acts faithfully on $X_K$.
Indeed, the effective model is a finite flat group scheme~$G'$
by theorem~\ref{representability_group_proper_case}. Since $G$
is reductive by the assumption on its order, we get $N=1$.
We prove a refinement of this result in
proposition~\ref{various_properties} below. There, we also give
other properties of the effective model of a finite
group scheme, especially in the case where the action is
{\em admissible}, which means that $X$ can be covered by $G$-stable
open affine subschemes. In this case, there exist quotient schemes
$X/G$ and $X/G'$, and we want to compare them.

\begin{prop} \label{various_properties}
Let $X$ be an $R$-scheme satisfying the assumptions of
theorem~\ref{representability_group} or of
theorem~\ref{representability_group_proper_case}.
Let $G$ be a finite flat $R$-group scheme acting on $X$ and
let $G'$ be its effective model. Then~:
\begin{trivlist}
\itemn{i} Let $W$ be a closed or an open subscheme of $X$. If
$W$ is $G$-stable, then it is $G'$-stable. In particular, if $G$
acts admissibly, then $G'$ also acts admissibly.
\itemn{ii} The effective model of a finite flat subgroup $H\subset G$,
for the restricted action on $X$, is the schematic image of $H$
in $G'$. If $H$ is normal in $G$, then $H'$ is normal in $G'$.
\itemn{iii} Assume that $G$ is {\'e}tale and let $p=\ch(k)$. Let
$N\lhd G$ be the (unique) subgroup of $G$ such that $N_k$ is the kernel
of the action on $X_k$. Then, the effective model of $N$ is a connected
$p$-group.
\end{trivlist}
In the sequel, we assume that $G$ acts admissibly on $X$.
\begin{trivlist}
\itemn{iv} The identity of $X$ induces an isomorphism $X/G\simeq X/G'$.
\itemn{v} Assume that there is an open subset $U\subset X$ which is
universally schematically dense, such that $G'$ acts
freely on $U$. Then for any closed normal subgroup $H\lhd G$, the effective
model of $G/H$ acting on $X/H$ is $G'/H'$.
\itemn{vi} Under assumptions {\rm (iii)} and {\rm (v)}, the group
$G'$ has a connected-{\'e}tale sequence
$$1\to N' \to G'\to G/N \to 1 \ .$$
\end{trivlist}
\end{prop}

\begin{proo}
(i) If $W$ is a closed subscheme of $X$, then it follows from the
general remarks of subsection~\ref{def_and_example} that the morphism
$G\times W\to W$ extends to a morphism $G'\times W\to W$. Now assume
that $W$ is open. It is enough to prove that the underlying set of
$W$ is stable under~$G'$. Let $w\in W$ be a point and let $\Omega$
be its orbit, by which we mean the schematic image of
$G\times\Spec(k(w))$ in $X$. This is a closed subscheme of $X$,
hence $G'$-stable. Since $\Omega\subset W$, it follows that $W$
is $G'$-stable.

\noindent (ii) This is clear.

\noindent (iii)
Since the composition $N_k\to N'_k\hookrightarrow \Aut_k(X_k)$ is trivial
as a morphism of sheaves, the morphism $N_k\to N'_k$ also is. Moreover,
$N\to N'$ is dominant and closed hence surjective. Hence $N'_k$ is
infinitesimal so $N'$ is a $p$-group. Let us show that it is connected. We may
and do assume that $R$ is henselian. Then $N'$ has a connected-{\'e}tale
sequence whose {\'e}tale quotient we denote by $N'_\et$. The composition
$t:N\to N'\to N'_\et$ is trivial on the special fibre. Moreover, $t$ is
determined by its restriction to the special fibre because it is a morphism
between {\'e}tale schemes. So it is globally trivial. As $t$ is dominant we
get $N'_\et=1$ thus $N'$ is connected.

\noindent (iv) The quotient $X\to X/G$ is described, locally on a $G$-stable
open affine $U=\Spec(A)$, by the invariant ring $A^G=\{\,a\in
A,\,\mu_G(a)=1\otimes a\,\}$ where $\mu_G\colon A\to RG\otimes A$ is the
coaction. Now $\mu_G$ factors through the coaction $\mu_{G'}$ corresponding to
the action of $G'$ :
$$
A\to RG'\otimes A\hookrightarrow RG\otimes A
$$
Therefore, $A^{G'}=\{\,a\in A,\,\mu_{G'}(a)=1\otimes a\,\}=A^G$. The result
follows.

\noindent (v) Clearly $H$ acts admissibly, and $X/H\simeq X/H'$
by (ii). We just have to show that $G'/H'$ acts faithfully on
$X/H'$. This is true since $G'/H'$ acts freely on the image of
$U$ in $X/H'$, by the assumptions on $U$.

\noindent (vi) Apply (v) to $H=N$.
\end{proo}

In~\ref{examples_torsors} and \ref{eff_model_of_quotient} below,
we will give an example where the effective model $G'$ does not
act freely on some schematically dense open subscheme, and the
claim in (v) does not hold.

\subsection{Schematic images for formal schemes}

The same methods as in subsection~\ref{main_result} yield analogous
representability results in the category of formal schemes locally of
finite type. Since the proofs are completely similar, we will simply
indicate how the objects are defined and then state the results. In
this subsection, the discrete valuation ring $(R,K,k,\pi)$ is complete
and we write $R_n:=R/\pi^n$. With a slight abuse, we use the notation
$i_n$ for both closed immersions $\Spec(R_n)\into\Spec(R_{n+1})$ and
$\Spec(R_n)\into \Spec(R)$, since confusions are not likely to arise.

\begin{noth} {\bf Formal sheaves.}
We first recall some notations and definitions. By a {\em presheaf}
over $R$ we mean a contravariant functor from the category of
$R$-schemes to the category of sets. As usual, we have the notion
of a {\em group presheaf} and most of what will be said hereafter
is valid for group presheaves. Schemes over~$R$ are identified with
their functor of points and hence can be viewed as presheaves.
Presheaves over $R$ form a category denoted $\PSh/R$. Of course
what we just said works for any base ring.

Let $i_n^*:\PSh/R_{n+1}\to\PSh/R_n$ be the pullback defined
by $i_n^*F=F\times_{\Spec(R_{n+1})}\Spec(R_n)$.
An {\em fppf formal sheaf over $R$} is a functor from the category
of formal $R$-schemes to the category of sets satisfying the sheaf
condition for fppf coverings. It may be identified with a direct
system of fppf sheaves over $R_n$, i.e. a sequence $(F_n)$ such
that $F_n=i_n^*F_{n+1}$ for all $n\ge 1$. Precisely, the identification
goes as follows~: to a formal sheaf $F$, we associate the direct
system $F_n=i_n^*F$. To a direct system $(F_n)$ of fppf sheaves
over $R_n$, we associate the functor $F=\varinjlim F_n$ defined
by $F(X)=\varprojlim F_n(X_n)$ where $X=(X_n)$.
These mappings are inverse to each other. We say that $F$ is locally
of finite presentation (or locally of finite type, since $R$ is
noetherian) if each $F_n$ is locally of finite presentation,
i.e. satisfies the usual condition of commutation with filtering
direct limits of rings (\cite{EGA}~IV.8.14.2).
\end{noth}

\begin{noth} {\bf Formal groups.}
Given formal $R$-schemes of finite type $X$ and $Y$, we have
two important examples of formal sheaves locally of finite type~: the
{\em homomorphism sheaf} $\Hom_R(X,Y)=\varinjlim
\Hom_{R_n}(X_n,Y_n)$ and the {\em automorphism sheaf}
$\Aut_R(X)=\varinjlim \Aut_{R_n}(X_n)$.

Let $G$ be a flat formal scheme in groups of finite type and $X$ a flat
separated formal scheme of finite type over~$R$. An action of
$G$ on $X$ is given by a morphism of formal schemes $G\times X\to X$
(satisfying the usual axioms) or equivalently by a morphism of formal
sheaves in groups $G\to \Aut_R(X)$. The kernel $N$ of the
action is defined as usual. As in
lemma~\ref{lemma_representability_of_kernel}, one shows
that $N$ is representable by a closed formal subscheme of $G$. As in
lemma~\ref{proper_gp_good_action}, one shows that there exists
a finite $R$-flat formal closed subscheme $Z\subset X$ such that
the induced action $G\times Z\to X$ has kernel equal
to~$G\times N$ (here the kernel is understood as a subobject
of $G\times G$, see subsection~\ref{kernels}). An action is
faithful if and only if $N=1$, and one can also define faithfulness
by requiring that no nontrivial $R$-flat closed subscheme of $G$
acts trivially on $X$.
\end{noth}

\begin{noth} {\bf Schematic images.}
Let $\Rig/K$ denote the category of quasi-compact,
quasi-separated rigid analytic $K$-spaces. As we recalled, Raynaud's
point of view gives an equivalence between $\Rig/K$ and the
category of flat formal $R$-schemes of finite type localised
by admissible formal blowing-ups. Using the existence of flat models
for flat morphisms of rigid spaces (see~\cite{BL2}), one can set up a
satisfactory theory of fppf descent in $\Rig/K$. It is not our
intention to provide the details of such a theory, as there
are more qualified experts to do this. We quote these facts
without further justification~; they give a meaning to what
an fppf sheaf on $\Rig/K$ is.

Recall that a {\em model} of a rigid $K$-space $X_K$ is a pair
$(X,i)$ where $X$ is a flat formal scheme of finite type and~$i$ is an
isomorphism between $X_\rig$ and $X_K$. A map between models
$(X_1,i_1)$ and $(X_2,i_2)$ is a morphism of formal schemes
$X_1\to X_2$ compatible with the given isomorphisms $i_1,i_2$.
We define the {\em generic fibre} $F_\rig$ of an fppf formal
sheaf locally of finite type $F$ to be the fppf sheaf on
$\Rig/K$ defined as follows. For any quasi-compact,
quasi-separated rigid analytic space $X_K$, we set~:
$$
F_\rig(X_K)=\varinjlim_{X_\rig=X_K} F(X)
$$
where the limit is taken with respect to all models $X$ of $X_K$.
If $F$ is representable by a formal scheme locally of finite type,
this definition coincides with the definition of the generic fibre of
a formal scheme by \cite{dJ}, proposition~7.1.7. Then the definitions
of the schematic closure of a subsheaf $G$ of the generic fibre
$F_\rig$, schematic image and related notions are the obvious
extensions of the definitions in subsection~\ref{def_and_example}. We
can now state our results for formal schemes.
\end{noth}

\begin{theo} \label{representability_formal_group}
Let $X$ be an affine flat formal $R$-scheme of finite type.
Let $G$ be a flat formal $R$-scheme in groups of finite type
acting on $X$, faithfully on the generic fibre. Then the
effective model $G'$ of the action is representable by a
flat formal $R$-scheme in groups which is not necessarily
of finite type. If $G$ is quasi-compact, or
affine, or finite, then $G'$ has the same property.
\end{theo}

\begin{theo} \label{representability_formal_group_proper_case}
Let $X$ be a flat, separated formal $R$-scheme of finite
type. Let $G$ be a proper flat formal $R$-scheme in groups
acting on $X$, faithfully on the generic fibre. Let $N$
denote the kernel of the action and assume that $N_k$ is
finite. Then the effective model $G'$ is representable by
a proper flat formal $R$-group scheme.
\end{theo}

\section{Examples} \label{Examples}

\subsection{Schematic closure of a $K$-group scheme}

When it is representable, it is clear that the schematic
image $G'$ depends only on
the generic fibre of $G$. One may start from an action of a
finite $K$-group scheme $G_K$ and wonder if its schematic closure
in $\Aut_R(X)$ is representable by a finite flat $R$-scheme.
This is not true in general, simply because the action of $G$
may fail to extend to the special fibre. For an example of  this,
consider the ring of power series $R=k[[\lambda]]$ over a field of
characteristic~$0$. Consider the projective completion of
the affine $R$-curve with equation $y^2=x(x-1)(x-\lambda)$, and let
$E/R$ be the complement of the unique singular point of the special
fibre. Thus $E_K$ is the Legendre elliptic curve over $K$.
The $2$-torsion $E_K[2]$ is rational and contains in particular
the point $A=(0,0)$ generating a group of translations
$G_K\simeq(\zmod{2})_K$. This point has singular reduction, and it is
easy to see that the image of the nontrivial point of $G_K$
under $G_K\to \Aut_R(E)$ is a closed point. Therefore, the
schematic closure is the group obtained by glueing $G_K$ and
the unit section $1_R$~; it is not finite over $R$.

\subsection{Two effective models of $\zmod{p^2}$}
\label{examples_torsors}

The end of the paper is devoted to the computation of schematic
images for the group $\zmod{p^2}$. The degeneration of torsors
under $\zmod{p}$ is well understood~; one observes the exceptional
feature that the effective model tends to act freely on an
$R$-universally dense open set. Recently, Sa{\"\i}di studied degenerations
of torsors under $\zmod{p^2}$ in equal characteristics~\cite{Sa}. He
computed equations for such degenerations~; they inherit an action of
$\zmod{p^2}$. We will compute the effective model in two cases~:
one case where one gets a torsor structure, and one where this
fails to happen. In the case of mixed characteristics,
similar examples have been given by Tossici in his Ph.D. thesis
using the Kummer-to-Artin-Schreier isogeny of Sekiguchi and Suwa
in degree $p^2$ (see \cite{To1}, \cite{To2}).

We let $(R,K,k,t)$ be a complete discrete valuation ring with equal
characteristics $p>0$, so $R\simeq k[[t]]$. Under this assumption, torsors
under $\zmod{p^2}$ are described by Witt theory.

\begin{noth}
{\bf Classical Witt theory.}
First we briefly recall the notations of Witt theory in degree $p^2$
(see \cite{DG}, chap. V). The group scheme of Witt vectors of length 2
over $R$ has underlying scheme
$W_{2,R}=\Spec(R[u_1,u_2])\simeq \bA^2_R$ with multiplication law
$$
(u_1,u_2)+(v_1,v_2)=\big(u_1+v_1,u_2+v_2+\sum_{k=1}^{p-1} \bin{p}{k}
\, u_1^kv_1^{p-k}\big) \ .
$$
Here we put once for all $\bin{p}{k}:=\frac{1}{p}{p\choose k}$ where
${p\choose k}$ is the binomial coefficient.
The Frobenius morphism of $W_2$ is denoted by $F(u_1,u_2)=(u_1^p,u_2^p)$.
Put $\phi:=F-\Id$. From the exact sequence
$$
0\to (\zmod{p^2})_R\to W_{2,R}\stackrel{\phi}{\longrightarrow} W_{2,R}\to 0
$$
it follows that any {\'e}tale torsor $f\colon \Spec(B)\to \Spec(A)$
under $(\zmod{p^2})_R$ is given by an equation
$$
F(X_1,X_2)-(X_1,X_2)=(a_1,a_2)
$$
where $(a_1,a_2)\in W_2(A)$ is a Witt vector and the substraction
is that of Witt vectors. Furthermore, $(a_1,a_2)$ is well-defined up
to addition of elements of the form $F(c_1,c_2)-(c_1,c_2)$. Note that
$$
F(X_1,X_2)-(X_1,X_2)=\big(X_1^p-X_1,X_2^p-X_2+
\sum_{k=1}^{p-1} \bin{p}{k} \, (X_1)^{pk}(-X_1)^{p-k}\big) \ .
$$
We emphasize that the Hopf algebra of $(\zmod{p^2})_R$ is
$$
R[\zmod{p^2}]=\frac{R[u_1,u_2]}{(u_1^p-u_1,u_2^p-u_2)}
$$
with comultiplication that of $W_2$.
\end{noth}

\begin{noth} \label{twistedforms}
{\bf Twisted forms of $W_2$.}
Let $\lambda,\mu,\nu$ be elements of $R$. We define a ''twisted''
group $W_2^\lambda$ as the group with underlying scheme
$\Spec(R[u_1,u_2])$ and multiplication law given by
$$
(u_1,u_2)+(v_1,v_2)=\bigg(u_1+v_1\, ,\, u_2+v_2+\lambda\,
\sum_{k=1}^{p-1} \bin{p}{k} \, u_1^kv_1^{p-k}\bigg) \ .
$$
We have the following analogues of the scalar multiplication and the
Frobenius of $W_2$ :
$$
\begin{array}{rl}
I_{\lambda,\mu}^{\nu} \, \colon \; W_2^\lambda & \longrightarrow W_2^{\lambda\mu}
\medskip \\
(u_1,u_2) & \longmapsto (\nu u_1,\mu\nu^p u_2) \\
\end{array}
$$
and
$$
\begin{array}{rl}
F_{\lambda} \, \colon \; W_2^\lambda & \longrightarrow W_2^{\lambda^p}
\medskip \\
(u_1,u_2) & \longmapsto (u_1^p,u_2^p) \ . \\
\end{array}
$$
In case $\mu=\lambda^{p-1}$ we define an isogeny
$$
\phi_{\lambda,\nu}:=F_{\lambda}-I_{\lambda,\lambda^{p-1}}^{\nu}\colon
W_2^\lambda \to W_2^{\lambda^p} \ .
$$
We have
$$
\phi_{\lambda,\nu}(u_1,u_2) =
\bigg(u_1^p-\nu u_1\, ,\, u_2^p-\nu^p\lambda^{p-1}u_2
+\lambda^p \sum_{k=1}^{p-1} \bin{p}{k} \, u_1^{pk}(-\nu u_1)^{p-k} \bigg) \ .
$$
The kernel $\cK_{\lambda,\nu}:=\ker(\phi_{\lambda,\nu})$
is a finite flat group of rank $p^2$. If $p>2$ its Hopf algebra is
$$
R[\cK_{\lambda,\nu}]=\frac{R[u_1,u_2]}{(u_1^p-\nu u_1,u_2^p-\nu^p\lambda^{p-1}u_2)} \ .
$$
\end{noth}


We now come to the examples. They arise from the following
situation. Denote by $G=\zmod{p^2}$ the constant group, and by
$Y=\bA^1_R=\Spec(R[w])$ the affine line
over $R$. Let $m_1,m_2\in \bZ$ be integers. Let
$f_K\colon X_K\to Y_K$ be the $(\zmod{p^2})_K$-torsor
over $Y_K=\bA^1_K$ given by the equations :
$$
\left\{
\begin{array}{lcl}
T_1^p-T_1 & = & t^{m_1}w \\
T_2^p-T_2 & = & t^{m_2}w - {\displaystyle \sum_{k=1}^{p-1}} \bin{p}{k} \, (T_1)^{pk}(-T_1)^{p-k} \\
\end{array}
\right.
$$
Depending on the values of the conductors
$m_1$, $m_2$ this gives rise to different group degenerations.

\begin{noth} 
{\bf First example.}
Assume $m_1=0$ and $m_2=-p$. Then after the change of variables $Z_1=T_1$,
$Z_2=t T_2$ the map $f_K$ extends to a cover $X\to Y$ with equations
$$
\left\{
\begin{array}{lcl}
Z_1^p-Z_1 & = & w \\
Z_2^p-t^{(p-1)}Z_2 & = & w - t^p{\displaystyle \sum_{k=1}^{p-1}} \bin{p}{k} \, (Z_1)^{pk}(-Z_1)^{p-k} \\
\end{array}
\right.
$$
The scheme $X$ is a smooth affine $R$-curve.
It is quickly seen that the action of $\zmod{p^2}$ extends to~$X$.
As is obvious from the expression of the isogeny $\phi_{\lambda,\nu}$
(see \ref{twistedforms}), the map $X\to Y$ is a torsor under
$\cK_{\lambda,\nu}$ for $\lambda=t$ and $\nu=1$. Thus, the effective
model is $G'=\cK_{t,1}$.
\end{noth}

\begin{noth} \label{second_example}
{\bf Second example.}
Assume $m_1=-p^2n_1<0$ and $m_2=0$. Put $\tilde{m}_1=n_1(p(p-1)+1)$.
Then after the change of variables $Z_1=t^{pn_1} T_1$ and
$Z_2=t^{\tilde{m}_1}T_2$ the map $f_K$ extends to a cover $X\to Y$
with equations
$$
\left\{
\begin{array}{lcl}
Z_1^p-t^{(p-1)pn_1}Z_1 & = & w \\
Z_2^p-t^{(p-1)\tilde{m}_1}Z_2 & = & t^{p\tilde{m}_1}w
- {\displaystyle \sum_{k=1}^{p-1}} \bin{p}{k} \, t^{pn_1(p-1)(p-1-k)}(Z_1)^{pk}(-Z_1)^{p-k} \ . \\
\end{array}
\right.
$$
The scheme $X$ is a flat $R$-curve with geometrically integral
cuspidal special fibre. The action of $\zmod{p^2}$ extends to
this model as follows : for $(u_1,u_2)$ a point of
$G_R=(\zmod{p^2})_R$,
$$
(u_1,u_2).(Z_1,Z_2)=
\bigg(Z_1+t^{pn_1}u_1\; ,\; Z_2+t^{\tilde{m}_1}u_2+
{\displaystyle \sum_{k=1}^{p-1}}
\bin{p}{k} \, t^{n_1(p(p-1)+1-pk)}(Z_1)^k(u_1)^{p-k}\bigg) \ .
$$
In order to find out the effective model $G'$ we look at the subalgebra
of $RG$ generated by $v_1=t^{n_1}u_1$ and $v_2=t^{\tilde{m}_1}u_2$ :
$$
RG':=R[v_1,v_2]\subset RG \ .
$$
One computes that $RG'$ inherits a comultiplication from $RG$ :
$$
(v_1,v_2)+(w_1,w_2)=\bigg(v_1+w_1\; ,\; v_2+w_2+
{\displaystyle \sum_{k=1}^{p-1}}
\bin{p}{k} \, t^{n_1(p-1)^2} v_1^kw_1^{p-k}\bigg) \ .
$$
Thus if $p>2$ we recognize $G'\simeq \cK_{\lambda,\nu}$ for
$\lambda=t^{n_1(p-1)^2}$ and $\nu=t^{n_1(p-1)}$. The action
of $G$ on $X$ extends to an action of $G'$ given by
$$
(v_1,v_2).(Z_1,Z_2)=
\bigg(Z_1+t^{(p-1)n_1}v_1\; ,\; Z_2+v_2+
{\displaystyle \sum_{k=1}^{p-1}}
\bin{p}{k} \, t^{n_1(p-1)(p-1-k)}Z_1^kv_1^{p-k}\bigg) \ .
$$
Here $X\to Y$ is not a torsor under $G'$. Indeed, on the special fibre
we have $G'_k=(\alpha_p)^2$ and the action on $X_k$ is
$$
(v_1,v_2).(Z_1,Z_2)=
\bigg(Z_1\; ,\; Z_2+v_2+v_1Z_1^{p-1}\bigg) \ .
$$
This action is faithful as required, but any point
$(z_1,z_2)\in X_k$ has a stabilizer of order $p$ which is the
subgroup of $G'_k$ defined by the equation $v_2+v_1 z_1^{p-1}=0$.
\end{noth}

\subsection{Effective model of a quotient} \label{eff_model_of_quotient}

We finish with a counter-example to point (v) in
proposition~\ref{various_properties}. For $\nu\in R$ we introduce
the group scheme $\cM_\nu$ which is the kernel of the isogeny
$\psi_\nu:\bG_{a,R}\to\bG_{a,R}$ defined by $\psi_\nu(x)=x^p-\nu x$
(see~\cite{Ma}, \S~3.2). This is a finite flat group scheme of order $p$.

We continue with the example in~\ref{second_example}. Thus
$G=(\zmod{p^2})_R$ and $G'\simeq \cK_{\lambda,\nu}$ where
$\lambda=t^{n_1(p-1)^2}$ and $\nu=t^{n_1(p-1)}$. Let
$H=(\zmod{p})_R\subset G$ and let $H'\subset G'$ be its image.
We have
$$
H'=\Spec\left(\frac{R[v_2]}{(v_2^p-\nu^p\lambda^{p-1}v_2)}\right)
\simeq \cM_{\nu^p\lambda^{p-1}}
$$
and
$$
G'/H'=\Spec\left(\frac{R[v_1]}{(v_1^p-\nu v_1)}\right)
\simeq \cM_{\nu} \ .
$$
The quotient scheme $X/H\simeq X/H'$ is the cover of $Y$ given
by the equation $Z_1^p-t^{(p-1)pn_1}Z_1=w$ i.e. $Z_1^p-\nu^pZ_1=w$.
It has and action of $G'/H'$ given by
$$
v_1.Z_1=Z_1+\nu v_1 \ .
$$
This action is not faithful on the special fibre. It is visible
that the effective model of $G'/H'$, or equivalently of $G/H$,
acting on $X/H'$ is the group whose Hopf algebra is equal to the
subalgebra of $R[G'/H']$ generated by $s_1=\nu v_1$. Therefore
$(G/H)'\simeq \cM_{\nu^p}$ and the map $G'/H'\to (G/H)'=(G'/H')'$
is not an isomorphism. We see that the effective model of the
quotient is not the quotient of the effective models.

\end{document}